\documentclass[12pt]{article}
\usepackage{latexsym}
\usepackage{amsmath}
\usepackage{amssymb}
\newcommand{\hs}{\hspace{0.5cm}}
\newcommand{\insp}{\hspace*{1cm}}

\def\ra{\rightarrow}

\newcommand{{\CB}}{\cal B}

\newcommand{{\CC}}{\cal C}

\newcommand{\od}{\overline{\Delta}}

\newcommand{\btr}{\blacktriangleright}
\newcommand{\btl}{\blacktriangleleft}

\def\ra{\mathop{\rightarrow}}

\def\picture#1 by #2 (#3){
\vbox to #2{
\hrule width #1 height 0pt depth 0pt
\vfill
\special{picture #3}}}

\begin{document}
\begin{center}
{\Large \bf The Drinfel'd double for group-cograded multiplier Hopf algebras}\\
\ \\
Lydia Delvaux\\
Department of Mathematics, L.U.C., Universiteitslaan, B-3590 Diepenbeek (Belgium)\\
E-mail: lydia.delvaux@luc.ac.be\\
\ \\
Alfons Van Daele\\
Department of Mathematics, K.U.Leuven, Celestijnenlaan 200B, B-3001 Heverlee (Belgium)\\
E-mail: alfons.vandaele@wis.kuleuven.ac.be
\end{center}
\ \\
\ \\

\begin{abstract}
\noindent Let $G$ be any group and let $K(G)$ denote the
multiplier Hopf algebra of complex functions with finite support
in $G$.  The product in $K(G)$ is pointwise.  The comultiplication
on $K(G)$ is defined with values in the multiplier algebra $M(K(G)
\otimes K(G))$ by the formula $(\Delta(f)) (p,q) = f(pq)$ for all
$f \in K(G)$ and $p, q \in G$. In this paper we consider
multiplier Hopf algebras $B$ (over $\Bbb C$) such that there is an
embedding $I: K(G) \rightarrow M(B)$.  This embedding is a
non-degenerate algebra homomorphism which respects the
comultiplication and maps $K(G)$ into the center of $M(B)$. These
multiplier Hopf algebras are called {\it $G$-cograded multiplier
Hopf algebras.} They are a generalization of the Hopf
group-coalgebras as studied by Turaev and Virelizier.
\\
In this paper, we also consider an {\it admissible} action $\pi$
of the group $G$ on a $G$-cograded multiplier Hopf algebra $B$.
When $B$ is paired with a multiplier Hopf algebra $A$, we
construct the Drinfel'd double $D^\pi$ where the coproduct and the
product depend on the action $\pi$. We also treat the $^*$-algebra
case.
\\
If $\pi$ is the trivial action, we recover the usual Drinfel'd
double associated with the pair $\langle A, B \rangle$. On the
other hand, also the Drinfel'd double, as constructed by Zunino
for a finite-type Hopf group-coalgebra, is an example of the
construction above. In this case, the action is non-trivial but
related with the adjoint action of the group on itself. Now, the
double is again a $G$-cograded multiplier Hopf algebra.
\end{abstract}
\ \\
{\it Mathematics Subject Classification}: 16W30, 17B37.
\\
\
\\
April 2004 ({\it Version 1.0})

\section*{Introduction}
Let $A$ be an algebra over $\Bbb C$.  If $A$ has no unit, we
require that the product in $A$ is non-degenerate as a bilinear
map. The multiplier algebra $M(A)$ of $A$ is the largest algebra
with unit in which $A$ sits as a dense two-sided ideal.  If $A$
has a unit, then $M(A) = A$.
\\
Consider a group $G$ and let $A$ be the algebra of complex valued
functions with finite support in $G$, with pointwise product. This
algebra has no unit, except when $G$ is finite. The multiplier
algebra $M(A)$ is given by the algebra of all complex valued
functions on $G$.  We define a comultiplication $\Delta$ on $A$ by
$(\Delta (f)) (p,q) = f(pq)$ where $f \in A$ and $p, q \in G$.  If
$G$ is finite, then $\Delta$ maps $A$ into $A \otimes A$ and makes
$A$ into a Hopf algebra. However, when $G$ is not finite,
$\Delta(f) \in M(A \otimes A)$ for all $f \in A$. In this case
$(A,\Delta)$ is a multiplier Hopf algebra, as reviewed in Section
1.  In this paper, the multiplier Hopf algebra $A$ associated with
a group $G$ as above, is denoted as $(K(G),\Delta)$.
\\
Multiplier Hopf algebras are generalizations of Hopf algebras when
the underlying algebra is no longer assumed to have a unit.
Integrals on multiplier Hopf algebras are defined as in the Hopf
algebra case, see Section 1. If $A$ is a multiplier Hopf algebra
with integrals, the dual object can be defined within the same
category.  This duality generalizes the one for finite-dimensional
Hopf algebras, but applies to a much bigger class of (multiplier)
Hopf algebras.  In fact, the theory of multiplier Hopf algebras is
a theory that allows results which are not possible within the
framework of usual Hopf algebras. Furthermore, this theory is also
a good model for an analytical theory of locally compact quantum
groups. The link between these two theories are the multiplier
Hopf $^\ast$-algebras with positive integrals, see [VD2], [VD3]
and especially [K-VD].
\\
Let $G$ be any group. In this paper we deal with {\it $G$-cograded
multiplier Hopf algebras} in the sense of Definition 2.1 (in
Section 2 of this paper). Roughly speeking, a multiplier Hopf
algebra $B$ is $G$-cograded if there is a central, non-degenerate
embedding $I : K(G) \rightarrow M(B)$. Furthermore, this embedding
respects the comultiplication, i.e. $\Delta (I(f)) = (I \otimes I)
(\Delta (f))$ for all $f \in K(G)$. Remark that we give a meaning
to this equation by extending the homomorphism $\Delta$ from $B$
to $M(B)$ and $I \otimes I$ from $K(G) \otimes K(G)$ to $M(K(G)
\otimes K(G))$ using the fact that the homomorphisms are
non-degenerate and hence have unique extensions.  It is shown in
[A-D-VD] that a Hopf group-coalgebra, as introduced by Turaev in
[T], is a special case of a group-cograded multiplier Hopf
algebra.  Therefore, we can interpret the results of Turaev,
Virelizier and Zunino within the theory of cograded multiplier
Hopf algebras.  This threws a new light on their results.  More
precisely, a lot of the results for Hopf group-coalgebras follow
from the more general results for multiplier Hopf algebras.
Moreover, we can apply the techniques from the theory of
multiplier Hopf algebra in the study of Hopf group-coalgebras. We
refer to [A-D-VD] for details about this approach to Hopf
group-coalgebras.
\\
The main goal of this paper is to apply this point of view when
constructing the quantum double for such $G$-cograded multiplier
Hopf algebras. We recover and generalize the work of Zunino on
this subject as it is found in [Z].
\\ \ \\
Let us now summarize the {\it content} of this paper.
\\
In {\it Section 1}, we recall the definition of a multiplier Hopf
algebra and we review some results  which are used in this paper.
\\
In {\it Section 2}, we first recall the notion of a group-cograded
multiplier Hopf algebra as studied in [A-D-VD]. Then we  consider
an {\it admissible} action of the group $G$ on a $G$-cograded
multiplier Hopf algebra $B$. Using this action, we deform the
comultiplication of $B$. This gives rise to a new multiplier Hopf
algebra $\widetilde{B} = (B,\widetilde{\Delta})$ where the
underlying algebra structure of $B$ is unchanged but with a
different comultiplication.  If $B$ is regular, so is $\widetilde
B$. If it is a  multiplier Hopf $^\ast$-algebra and has positive
integrals, then the same is true for $\widetilde{B}$, see Theorem
2.11.
\\
In {\it Section 3}, we start with a multiplier Hopf algebra
pairing $\langle A, B\rangle$ where $B$ is a $G$-cograded
multiplier Hopf algebra.  Let $\pi$ be an admissible action of $G$
on $B$.  Then we construct a twisted tensor product multiplier
Hopf algebra (as reviewed in Section 1) of the multiplier Hopf
algebras $(A, \Delta^{cop})$ and $(B, \widetilde{\Delta})$.  This
twisted tensor product multiplier Hopf algebra is denoted as
$D^{\pi} = A^{cop} \bowtie \widetilde{B}$. The main results on
$D^\pi$ are given in Theorem 3.8 and Proposition 3.9. If $\pi$ is
the trivial action, we recover the usual Drinfel'd double of the
pair $\langle A,B\rangle$, as reviewed in Section 1.  The
Drinfel'd double, as constructed by Zunino for a finite-type Hopf
group-coalgebras in [Z] is considered in Example 3.14.
\\
By the constructions in Section 2 and Section 3, we have further
interesting, non-trivial examples of multiplier Hopf
($^\ast$-)algebras with (positive) integrals.
\\ \ \\
Let us finish this introduction by mentioning some {\it basic
references}.  For multiplier Hopf algebras, these are  [VD1] and
[VD-Z].  For multiplier Hopf algebras with integrals, we refer to
[VD2].  The Drinfel'd double of a pair of multiplier Hopf algebras
is studied extensively in [Dr-VD] and in [De-VD]. The Hopf
group-coalgebras are introduced in [T]. The new approach to Hopf
group-coalgebras is studied and generalized to the case of
multiplier Hopf algebras in [A-D-VD].
\\
\ \\
{\bf Acknowledgements}\\
We are indebted to A.T.\ Abt El-Hafez, A.S.\ Hegazi and M.\
Mansour who provided us with their preprint ``Multiplier Hopf
group-coalgebras'' [H-A-M] and in doing so, brought to our
attention the works of Turaev, Virelizier and Zunino.

\section{Preliminaries on multiplier Hopf ($^\ast$-)algebras}\label{preliminaries}

We begin this section with a short introduction to the theory of
multiplier Hopf algebras.
\\ \ \\
{\it Multiplier Hopf ($^\ast$-)algebras with (positive)
integrals}\\
\ \\
As mentioned already in the introduction, throughout this paper,
all algebras are algebras over the field $\Bbb C$ of complex
numbers.  They may or may not have units, but always should be
non-degenerate, i.e. the multiplication maps (viewed as bilinear
forms) are non-degenerate.  For an algebra  $A$, the multiplier
algebra $M(A)$ of $A$ is defined as the largest algebra with unit
in which $A$ is a dense ideal, i.e. $A$ has no (left and right)
annihilators in $M(A)$.
\\
Now, we recall the definition of a multiplier Hopf algebra (see
[VD1] for details). Consider the tensor product $A \otimes A$
which is again an algebra with a non-degenerate product.  The
embedding of $A \otimes A\hookrightarrow M(A \otimes A)$ factors
through $M(A) \otimes M(A)$ in an obvious way as follows: $A
\otimes A \hookrightarrow M(A) \otimes M(A) \hookrightarrow M(A
\otimes A)$. A comultiplication on $A$ is a homomorphism $\Delta :
A \rightarrow M(A \otimes A)$ such that $\Delta(a) (1 \otimes b)$
and $(a \otimes 1) \Delta(b)$ are elements of $A \otimes A$ for
all $a,b \in A$.  We require $\Delta$ to be coassociative in the
sense that
$$(a \otimes 1 \otimes 1) (\Delta \otimes \iota) (\Delta(b) (1\otimes
c)) = (\iota \otimes \Delta)((a \otimes 1)\Delta(b)) (1\otimes 1
\otimes c)$$
for all $a,b, c\in A$ (where $\iota$ denotes the identity map).
\\ \ \\
{\bf 1.1 Definition [VD1]} \hs
A pair $(A,\Delta)$ of an algebra $A$ with non-degenerate
pro\-duct and a comultiplication $\Delta$ on $A$ is called a
multiplier Hopf algebra if the linear maps $T_1, T_2 : A \otimes A
\rightarrow A \otimes A$, defined by
\begin{eqnarray*}
T_1 (a \otimes b) = \Delta(a) (1 \otimes b) \insp T_2 (a \otimes b)
= (a \otimes 1) \Delta (b)
\end{eqnarray*}
are bijections.\\
\ \\
The conditions in Definition 1.1 in fact imply that $\Delta$ is a
non-degenerate homomorphism. For the convenience of the reader, we
recall the notion of a {\it non-degenerate homomorphism}. It is a
homomorphism $\gamma:A\to M(B)$, where $A$ and $B$ are algebras
with a non-degenerate product, such that
$\gamma(A)B=B\gamma(A)=B$. So, e.g.\ every element $b\in B$ is a
sum of elements of the form $\gamma(a)b'$ with $a\in A$ and $b'\in
B$. An important property of such a non-degenerate homomorphism
$\gamma$ is that it has a unique extension, to a unital
homomorphism from $M(A)$ to $M(B)$. The extension is still denoted
by $\gamma$. See the appendix in [VD1]. The homomorphisms $\iota
\otimes \Delta$ and $\Delta \otimes \iota$ are also non-degenerate
and so have unique extensions to $M(A \otimes A)$ in a natural
way. The coassociativity as formulated above means nothing else
but $(\Delta \otimes \iota) \Delta = (\iota \otimes \Delta)\Delta$
in
$M(A \otimes A \otimes A)$.\\

The bijectivity of the two maps in Definition 1.1 is equivalent
with the existence of a counit $\varepsilon$ and an antipode $S$
satisfying (and defined) by
\begin{eqnarray*}
\begin{array}{ll}
(\varepsilon \otimes \iota) (\Delta(a)(1 \otimes b)) = ab &\insp
m((S \otimes \iota)(\Delta(a) (1 \otimes b)))
= \varepsilon(a) b\\
(\iota \otimes \varepsilon)((a\otimes 1) \Delta (b)) = ab &\insp
m((\iota \otimes S)((a \otimes 1) \Delta(b))) = \varepsilon(b) a
\end{array}
\end{eqnarray*}

where $\varepsilon : A \rightarrow \Bbb C$ is a homomorphism, $S :
A \rightarrow M(A)$ is an anti-homomorphism and $m$ is the
multiplication map, considered as a linear map from $A \otimes A$
to $A$ and extended to $A \otimes M(A)$ and $M(A) \otimes A$.
\\ A multiplier Hopf algebra is called {\it regular} if $(A,
\Delta^{cop})$ is (also) a multiplier Hopf algebra, where
$\Delta^{cop}$ denotes the co-opposite comultiplication defined as
$\Delta^{cop} = \sigma \circ \Delta$ with $\sigma$ the usual flip
map from $A\otimes A$ to itself (and extended to $M(A \otimes
A)$). In this case, we also have that $\Delta (a) (b \otimes 1)$
and $(1\otimes b) \Delta(a)$ are in $A \otimes A$ for all $a, b
\in A$. A multiplier Hopf algebra is regular if and only if the
antipode is a bijection from $A$ to $A$ (see [VD2, Proposition
2.9]). Any Hopf algebra is a multiplier Hopf algebra. Conversely,
a multiplier Hopf algebra with unit is a Hopf algebra.
\\
In [Dr-VD], the use of the Sweedler notation for regular
multiplier Hopf algebras has been introduced.  We will also use
this notation in this paper. We will e.g. write $\sum a_{(1)}
\otimes a_{(2)} b$ for $\Delta(a) (1 \otimes b)$ and $\sum
ab_{(1)} \otimes b_{(2)}$ for $(a \otimes 1) \Delta(b)$.
\\ \ \\
{\bf 1.2 Definition [VD1]} \hs If $A$ is a $^\ast$-algebra, we
require the comultiplication $\Delta$ to be also a
$^\ast$-homomorphism. Then, a multiplier Hopf $^\ast$-algebra is a
$^\ast$-algebra with a comultiplication, making it into a
multiplier Hopf algebra. For a multiplier Hopf $^\ast$-algebra,
regularity is automatic.
\\
Recall that a $^\ast$-algebra $A$ over $\Bbb C$ is an algebra with
an involution $a \mapsto a^\ast$.  An involution is a antilinear
map satisfying $a^{\ast\ast} = a$ and $(ab)^\ast = b^\ast a^\ast$
for all $a$ and $b$ in $A$. The multiplier algebra $M(A)$ is again
a $^*$-algebra.
\\ \ \\
{\bf 1.3 Example} \hs Let $G$ be any group and let $A$ be the
$^\ast$-algebra $K(G)$ of complex, finitely supported functions on
$G$.  In this case $M(A)$ consists of all complex functions on
$G$.  Moreover $A \otimes A$ can be naturally identified with
finitely supported complex functions on $G \times G$ so that $M(A
\otimes A)$ is the space of all complex functions on $G \times G$.
If we define $\Delta : A \rightarrow M(A \otimes A)$ by $(\Delta
(f)) (p,q) = f(pq)\ \mbox{for}\ f \in A\ \mbox{and}\ p,q \in G$,
we clearly get a $^\ast$-homomorphism. If $f,g \in A$, then $(p,q)
\mapsto f(pq) g(q)$ and $(p,q) \mapsto g(p) f(pq)$ have finite
support and so belong to $A \otimes A$. The coassociativity
condition on $\Delta$ is a consequence of the associativity of the
multiplication on $G$.  So $\Delta$ is a comultiplication. To
obtain that the pair $(A,\Delta)$ is a multiplier Hopf algebra in
the sense of Definition 1.1, we notice that the bijectivity of the
linear maps $T_1$ and $T_2$ follows from the fact that the maps
$(p,q) \mapsto (pq,q)$ and $(p,q) \mapsto (p,pq)$ are bijective
from $G\times G$ to itself (because $G$ is assumed to be a group).
\\
This is a very simple example.  Interesting examples of multiplier
Hopf algebras are found among the discrete quantum groups (i.e.\
the duals of compact quantum groups), see e.g. [VD2].
\\ \ \\
We now discuss the notion of an {\it integral} on a multiplier
Hopf ($^\ast$-)algebra.  It is like the integral on a Hopf
algebra. We will restrict to the case of a regular multiplier Hopf
algebra (in particular to a multiplier Hopf $^\ast$-algebra).  We
first observe the following.  If $f$ is a linear functional on a
multiplier Hopf algebra $A$, we can define for all $a \in A$, the
multipliers $(\iota \otimes f) \Delta(a)$ in $M(A)$ in the
following way
\begin{eqnarray*}
((\iota \otimes f)\Delta(a)) b = (\iota \otimes f)(\Delta(a) (b
\otimes 1)) \insp b((\iota \otimes f) \Delta(a)) = (\iota \otimes
f) ((b\otimes 1) \Delta (a))
\end{eqnarray*}
where $b \in A$.  This is well-defined as both $\Delta(a)
(b\otimes 1)$ and $(b\otimes 1) \Delta(a)$ are in $A \otimes A$
and we can apply $\iota \otimes f$ mapping $A \otimes A$ to $A
\otimes \Bbb C$ (which is naturally itendified with $A$ itself).
Similarly, we can define $(f \otimes \iota) \Delta(a)$ in $M(A)$.
Then, the
following definition makes sense.
\\ \ \\
{\bf 1.4 Definition [VD2]} \hs A linear functional $\varphi$ on
$A$ is called left invariant if \\ $(\iota \otimes \varphi)
\Delta(a) = \varphi(a) 1$ for all $a \in A$.  A left integral is a
non-zero left invariant functional on $A$.  Similarly, a non-zero
linear functional $\psi$ satisfying $(\psi \otimes \iota)
\Delta(a) = \psi (a) 1$ for all $a \in A$ is called a right integral.
\\ \ \\
In Example 1.3, a left integral is given by the formula
$\varphi(f) = \sum_{q \in G} f(q)$ (the sum is well-defined as
only finitely many entries are non-zero).  In this example the
left integral is also right invariant.  This however is no longer
true in general.  Multiplier Hopf ($^\ast$-)algebras with
integrals are studied intensively in [VD2].  There are various
data (and many relations among them) about left and right
integrals. We collect some important results of [VD2].
\\ \ \\
{\bf 1.5 Theorem [VD2]} \hs Let $(A, \Delta)$ be a multiplier Hopf
($^\ast$-)algebra with left integral $\varphi$.  Any other left
integral is a scalar multiple of $\varphi$.  There is also a right
integral $\psi$, unique up to a scalar.    The left integral is
faithful in the sense that when $a \in A$, then $a = 0$ if
$\varphi(ab) = 0$ for all $b$ or $\varphi(ba) = 0$ for all $b$.
Similarly, the right integral is faithful. There is an
automorphism $\sigma$ of $A$ such that $\varphi(ab) =
\varphi(b\sigma(a))$ for all $a, b \in A$. There is an invertible
multiplier $\delta$ in $M(A)$ such that $(\varphi \otimes \iota)
\Delta(a) = \varphi(a) \delta$ and $(\iota \otimes \psi) \Delta(a)
= \psi(a) \delta^{-1}$ for all $a \in A$.\insp $\blacksquare$
\\ \ \\
One of the main features of a multiplier Hopf algebra $A$ with
integrals is the existence of the dual multiplier Hopf algebra
$(\widehat{A}, \widehat{\Delta})$. It is constructed in the
following way.
\\ \ \\
{\bf 1.6 Definition [VD2]} \hs Let $(A, \Delta)$ be a multiplier
Hopf ($^\ast$-)algebra with left integral $\varphi$.  Denote by
$\widehat{A}$ the space of linear functionals on $A$ of the form
$x \mapsto \varphi(xa)$ where $a \in A$.  The product
(respectively the coproduct) of $\widehat{A}$ is dual to the
coproduct (respectively the product) of $A$.
\\
\ This dual object $(\widehat{A}, \widehat{\Delta})$ is again a
regular multiplier Hopf algebra with integrals. Moreover, the dual
of $(\widehat{A}, \widehat{\Delta})$ is canonically isomorphic
with the orginial multiplier Hopf algebra $(A,\Delta)$. For a
finite-dimensional Hopf algebra $A$, we notice that $\widehat{A}$
equals the usual dual Hopf algebra of $A$. It is possible to
generalize many aspects of harmonic analysis in this general
framework. One can define the Fourier transform, one can prove
Plancherel's formula, ... . For details, see [VD2].
\\ \ \\
In the case that $(A,\Delta)$ is a multiplier Hopf
$^\ast$-algebra, a left integral is called positive if
$\varphi(a^\ast a) \geq 0$ for all $a \in A$.  We mention that,
when there is a positive left integral, there is also a positive
right integral.  Multiplier Hopf $^\ast$-algebras with positive
integrals, give rise to dual multiplier Hopf $^\ast$-algebras with
positive integrals. In [K-VD] is explained how any multiplier Hopf
$^\ast$-algebra with positive integrals gives rise to a locally
compact quantum group as introduced and studied by Kustermans and
Vaes in [K-V].
\\ \ \\
{\it Pairing and Drinfel'd double of multiplier Hopf algebras}
\\ \ \\
We now recall how the Drinfel'd double is constructed from a pair
of multiplier Hopf algebras.
\\
Start with two regular multiplier Hopf algebras $(A,\Delta)$ and
$(B,\Delta)$ together with a non-degenerate bilinear map $\langle
\cdot, \cdot \rangle$ from $A \times B$ to $\Bbb C$. This is
called a pairing if certain conditions are fulfilled. The main
property is that the coproduct in $A$ is dual to the product in
$B$ and vice versa. There are however certain regularity
conditions, needed to give a correct meaning to this statement.
The investigation of these conditions is done in [Dr-VD]. We
recall some important aspects here.
\\ \ \\
{\bf 1.7 Definition} \hs  For $a \in A$ and $b \in B$, define $a
\btr b$, $b \btl a$, $b \btr a$ and $a \btl b$ as multipliers in
the following way. Take $a' \in A$, $b' \in B$ and define
\begin{eqnarray*}
\begin{array}{lll}
(b \btr a) a' &= &\sum \langle a_{(2)}, b\rangle a_{(1)} a'\\
(a \btl b) a' &= &\sum \langle a_{(1)}, b\rangle a_{(2)} a'
\end{array}
\insp
\begin{array}{lll}
(a \btr b) b' &= &\sum \langle a, b_{(2)}\rangle b_{(1)} b'\\
(b \btl a) b' &= &\sum \langle a, b_{(1)}\rangle b_{(2)} b'.
\end{array}
\end{eqnarray*}
The regularity conditions on the pairing say (among other things)
that the multipliers $b \btr a$ and $a\btl b$ in $M(A)$ (resp. $a
\btr b$ and $b \btl a$ in $M(B)$) actually belong to $A$ (resp.\ $B$).
\\ \ \\
Then it is possible to state that the product and the coproduct
are dual to each other:
\begin{eqnarray*}
\begin{array}{lll}
\langle a, bb'\rangle &= &\langle b' \btr a,b\rangle\\
&= &\langle a \btl b, b'\rangle
\end{array}
\insp
\begin{array}{lll}
\langle aa', b\rangle &= &\langle a, a' \btr b\rangle\\
&= &\langle a', b \btl a\rangle.
\end{array}
\end{eqnarray*}
There are four modules involved.  All these modules are unital. By
definition, e.g. $B$ is a left $A$-module for the action $A \btr
B$.  That $B$ is unital means that any element $b \in B$ is a
linear combination of elements of the form $a \btr b'$ with $a \in
A$ and $b' \in B$. A stronger result however is possible here,
coming from the existence of local units, see [D-VD-Z].  Take e.g.
$b \in B$. Then there are elements $\{a_1, a_2, \ldots, a_n\}$ in
$A$ and $\{b_1, b_2, \ldots, b_n\}$ in $B$ such that $b = \sum a_i
\btr b_i$.  Because of the existence of local units, there is an
$e \in A$ such that $ea_i = a_i$ for all $i$.  It follows easily
that $e \btr b =b$. So, we have that for all $b \in B$ there
exists an element $e \in A$ such that $b = e \btr b$.
\\
As an important
consequence of this last result, we can use the Sweedler notation
in the framework of dual pairs in the following sense. Take $a \in
A$ and $b \in B$, and e.g. the element $b \btr a = \sum \langle
a_{(2)}, b\rangle a_{(1)}$.  In the right hand side of this
equation, the element $a_{(2)}$ is covered by $b$ through the
pairing because $b = e\btr b$ for some $e \in A$ and therefore
$\sum \langle a_{(2)}, b\rangle a_{(1)} = \langle a_{(2)}, e \btr
b \rangle a_{(1)} = \sum \langle a_{(2)} e, b\rangle a_{(1)}$.
\\
We also mention that $\langle S(a), b\rangle = \langle a,
S(b)\rangle$ and as expected, $\langle a, 1\rangle = \varepsilon
(a)$ and $\langle 1, b\rangle = \varepsilon (b)$ where $a \in A$
and $b \in B$.  For these formulas, one has to extend the pairing
to $A \times M(B)$ and to $M(A) \times B$. This can be done in a
natural way using the fact that the four modules $A \btr B$, $B
\btr A$, $A \btl B$ and $B \btl A$ are unital.\\ If $A$ and $B$
are multiplier Hopf $^\ast$-algebras.  A pairing $\langle A,
B\rangle$ is called a multiplier Hopf $^\ast$-algebra pairing if
additionally $\langle a^\ast, b\rangle = \overline{\langle
a,S(b)^\ast\rangle}$ for all $a \in A$ and $b \in B$.
\\ \ \\
A pairing of two multiplier Hopf algebras is the natural setting
for the construction of the Drinfel'd double and it turns out that
the conditions on the pairing $\langle A, B \rangle$ are
sufficient to make this construction. This is done in a rigorous
way in the papers [Dr-VD], [D] and [De-VD]. We recall some
essential ideas.
\\
The main point  is that there is (as in the case of
finite-dimensional Hopf algebras) an invertible twist map.
\\ \ \\
{\bf 1.8 Definition} \hs For $a \in A$ and $b \in B$, we set
\begin{eqnarray*}
R(b \otimes a) = \sum(b_{(1)} \btr a \btl S^{-1} (b_{(3)}))
\otimes b_{(2)}.
\end{eqnarray*}
It is proven in [Dr-VD] that this map is well-defined and
bijective. Let $D = A \bowtie B$ denote the algebra with the
tensor product $A \otimes B$ as the underlying space and with the
twisted product given by the twist map $R$ as follows:
\begin{eqnarray*}
(a \bowtie b) (a' \bowtie b') = (m_A \otimes m_B) (\iota_A \otimes
R \otimes \iota_B) (a \otimes b \otimes a' \otimes b')
\end{eqnarray*}
with $a, a' \in A$ and $b, b' \in B$. The maps $A \rightarrow M(D)
: a \mapsto a \bowtie 1$ and $B \rightarrow M(D) : b \mapsto 1
\bowtie b$ are non-degenerate algebra embeddings. The embedding of
$A$ in $M(D)$ gives rise to the embedding of $A \otimes A$ in $M(D
\otimes D)$. Similarly, $B \otimes B$ can be embedded in $M(D
\otimes D)$. These embeddings can be extended to the multiplier
algebras. The comultiplication on $D$ can then be given by the
formula $\Delta_D (a \bowtie b) = \Delta^{cop} (a) \Delta (b)$.
\\
The main result is the following.
\\ \ \\
{\bf 1.9 Theorem} \hs With the notations and definitions above,
the pair $(D,\Delta_D)$ is a multiplier Hopf algebra, called the
quantum double (or Drinfel'd double). If $A$ and $B$ have
integrals, then $D$ has integrals too.
\\
More precisely, let $\varphi_A$ denote a left integral on $A$ and
let $\psi_B$ denote a right integral on $B$, then $\psi_D =
\varphi_A \otimes \psi_B$ is a right integral on $D$.
\\
If $\langle A,B \rangle$ is a multiplier Hopf $^\ast$-algebra
pairing of two multiplier Hopf $^\ast$-algebras, then $D = A^{cop}
\bowtie B$ is again a multiplier Hopf $^\ast$-algebra.  Suppose
that $\varphi_A$ (resp. $\psi_B$) is a positive left integral on
$A$ (resp. right integral on $B$).  In [De-VD] it is shown that
there esists a complex number $\rho$ such that $\rho(\varphi_A
\otimes \psi_B)$ is a positive right integral on $D = A^{cop}
\bowtie B$.
\\ \ \\
{\it Twisted tensor product construction of multiplier Hopf
($^\ast$-)algebras}
\\ \ \\
The Drinfel'd double construction, as reviewed above, is a special
case of a twisted tensor product of multiplier Hopf
($^\ast$-)algebras.  General twisted tensor products of multiplier
Hopf algebras are studied in [D].  Also here, we recall some
ideas.
\\ \ \\
{\bf 1.10 Assumptions} \hs Let $A$ and $B$ be two algebras and
suppose that there is given a bijective linear map $R: B\otimes A
\to A\otimes B$ such that
\begin{eqnarray*}
\begin{array}{l}
R(m_B\otimes \iota_A)=(\iota_A\otimes m_B)
(R\otimes \iota_B)(\iota_B \otimes R)\\
R(\iota_B\otimes m_A)=(m_A\otimes \iota_B)(\iota_A\otimes R)
(R\otimes \iota_A).
\end{array}
\end{eqnarray*}
Recall that $m_A$ denotes the product in $A$, considered as a
linear map $m_A:A\otimes A \to A$ and similarly for the product
$m_B$ on $B$.
\\
One can consider the twisted tensor product
algebra $A\bowtie B$ in the following way. As a vector space $A
\bowtie B$ is $A \otimes B$. The product in $A \bowtie B$ is
defined by
\begin{eqnarray*} (a\bowtie b)(a' \bowtie b')=(m_A \otimes
m_B)(\iota_A \otimes R\otimes \iota_B)(a\otimes b \otimes
a'\otimes b')
\end{eqnarray*} for $a, a' \in A$ and $b, b' \in B$. As before,
we use  $\iota$ to denote the identity map, in particular we have
$\iota_A$ and $\iota_B$ for the identity maps on $A$ and $B$
respectively.
\\
The above assumptions on R are necessary for the associativity of
the product in $A\bowtie B$. Because the the products in $A$ and
$B$ are assumed to be non-degenerate, one can prove that the
bijectivity of $R$ guaranties that the product in $A\bowtie B$ is
again non-degenerate.
\\ \ \\
{\bf 1.11 Remarks} For details, see [D].
\begin{itemize}
\item[(i)] By using the conditions on $R$ and the bijectivity,
one can prove that the product in $A\bowtie B$ is also given
by the following expressions \begin{eqnarray*}
\begin{array}{l}
(a\bowtie b)(a'\bowtie b')=((\iota_A\otimes m_B)\circ
R_{12}\circ(\iota_B\otimes m_A\otimes
\iota_B)\circ(R^{-1})_{12})(a\otimes b\otimes a'\otimes b')\\
(a\bowtie b)(a'\bowtie b')=((m_A\otimes \iota_B)\circ
R_{34}\circ(\iota_A\otimes m_B\otimes
\iota_A)\circ(R^{-1})_{34})(a\otimes b\otimes a'\otimes b')
\end{array} \end{eqnarray*}
for all $a, a'\in A$ and $b,b'\in B$. Recall that we use the {\it
leg-numbering notation} for the maps $R$. When we write e.g.\
$R_{12}$, we consider $R$ as acting on the first two factors in
the tensor product. Similarly $R_{34}$ is the map $R$ as acting on
the 3th and the 4th factor. These formulas will be used to justify
that the decompositions of the comultiplications in $A$ and in $B$
are well covered when we use the Sweedler notation.
\item[(ii)] The maps
\begin{eqnarray*}
A \rightarrow M(A \bowtie B): a \mapsto a \bowtie 1 \insp B
\rightarrow M(A \bowtie B): b \mapsto 1 \bowtie b
\end{eqnarray*}
are non-degenerate algebra embeddings and therefore extend, in a
natural way, to unital algebra embeddings from $M(A)$ and $M(B)$
respectively to $M(A \bowtie B)$. Therefore, we have the
non-degenerate algebra embeddings
\begin{eqnarray*}
\begin{array}{l} A \otimes A \rightarrow M((A \bowtie B) \otimes
(A \bowtie B)) : a \otimes a' \mapsto (a \bowtie 1) \otimes (a'
\bowtie 1)\\ B \otimes B \rightarrow M((A \bowtie B) \otimes (A
\bowtie B)) : b \otimes b' \mapsto (1 \bowtie b) \otimes (1'
\bowtie b').
\end{array} \end{eqnarray*}
Also these embeddings
extend to the multiplier algebras in a natural way.
\item[(iii)]
Let $A$ and $B$ be *-algebras with non-degenerate products.
Suppose that the twist map $R:B\otimes A\to A\otimes B$ is
bijective and satisfies the above conditions. If furthermore
$(R\circ \ast_{B}\otimes\ast_{A}\circ\sigma)(R\circ
\ast_{B}\otimes\ast_{A}\circ\sigma)=\iota_A\otimes \iota_B$, then
there is a $^\ast$-operation on $A\bowtie B$ given as follows:
$(a\bowtie b)^*=R(b^*\otimes a^*)$ for all $a\in A, b\in B$. Now
the embeddings in Remark (ii) become $^\ast$-embeddings. \\
\end{itemize}
The comultiplications on $A$ and $B$ can be used to define the
comultipication on $A \bowtie B$ as usual:
\\ \ \\
{\bf 1.12 Definition [D]} \hs Let $A$ and $B$ be multiplier Hopf
algebras.  Let $R : B \otimes A \ra A \otimes B$ be a bijective
map satisfying the Assumptions 1.10. For $a \in A$ and $b \in B$,
define $$\overline{\Delta}(a \bowtie b) = \Delta(a) \Delta(b) \in
M((A \bowtie B) \otimes (A \bowtie B)).$$ In the next theorem, we
formulate sufficient conditions for $(A \bowtie B,
\overline{\Delta})$ to be a regular multiplier Hopf algebra.
\\ \ \\
{\bf 1.13 Theorem [D]} \hs Let $A$ and $B$ be multiplier Hopf
algebras with a bijective twist map $R$, satisfying the following
conditions
\begin{itemize}
\item[(1)] $R(m_B \otimes \iota_A) = (\iota_A \otimes m_B)(R \otimes \iota_B)(\iota_B \otimes R)$\\
$R(\iota_B \otimes m_A) = (m_A \otimes \iota_B)(\iota_A \otimes
R)(R \otimes \iota_A)$
\item[(2)] $\od (R(b \otimes a))= \Delta(b) \Delta(a)$ in
$M((A \bowtie B) \otimes (A \bowtie B))$ for all $a\in A$ and $b
\in B$.
\end{itemize}
Then $(A \bowtie B, \od, \overline{\varepsilon}, \overline{S})$ is
a regular multiplier Hopf algebra with, $\overline{\varepsilon}$
and $\overline{S}$ given as $\overline{\varepsilon} =
\varepsilon_A \otimes \varepsilon_B\ \mbox{and}\ \overline{S} =
R\circ (S_B \otimes S_A) \circ \sigma$.\\ Let $A$ and $B$ be
multiplier Hopf $^\ast$-algebras and $(R \circ (\ast_B \otimes
\ast_A) \circ \sigma) (R \circ (\ast_B \otimes \ast_A) \circ
\sigma) = \iota_A \otimes \iota_B$. Then $(A \bowtie B,
\overline{\Delta})$ is made into a multiplier Hopf
$^\ast$-algebra, if the $^\ast$-operation is defined by $(a
\bowtie b)^\ast = R(b^\ast \otimes a^\ast)$. \hfill$\blacksquare$

\section{Actions of $G$ on a $G$-cograded multiplier Hopf algebra}
Throughout this section, $G$ is an arbitrary group.  Let $K(G)$ be
the multiplier Hopf algebra of the complex valued functions with
finite support in $G$, see Example 1.3.
\\ \ \\
{\it Group-cograded multiplier Hopf algebras}
\\ \ \\
Recall the following definition given in [A-D-VD].

{\bf 2.1 Definition} \hs A {\it $G$-cograded multiplier Hopf
($^\ast$-)algebra} is multiplier Hopf ($^\ast$-)algebra $B$ so
that the following hold:
\begin{itemize}
\item[(1)] $B = \bigoplus\limits_{p \in G} B_p$ with $\{B_p\}_{p
\in G}$ a family of ($^\ast$-)subalgebras such that $B_p B_q = 0$
if $p \neq q$,
\item[(2)] $\Delta(B_{pq})(1 \otimes B_q) = B_p \otimes B_q$
 and $(B_p \otimes 1)\Delta (B_{pq}) = B_p \otimes B_q$
for all $p, q \in G$.
\end{itemize}
We say that the comultiplication of $B$ is $G$-graded.
\\ \ \\
{\bf 2.2 Proposition}\hs The data of a Hopf group-coalgebra, as
introduced by Turaev in [T], give an example of a cograded
multiplier Hopf
algebra.\\
For the proof of this propisition, we refer to [A-D-VD].\\
\ \\
In order to characterize the coalgebra structure of a general
$G$-cograded multiplier Hopf algebra, we first need the following
lemma.\\
\ \\
{\bf 2.3 Lemma} \hs Let $A$ and $B$ be multiplier Hopf algebras.
Let $f : A \rightarrow M(B)$ be a non-degenerate algebra
homomorphism which respects the comultiplication in the sense that
$\Delta_B \circ f = (f \otimes f) \circ \Delta_A$.  Then $f$
preserves the unit, the counit and the antipode in the following
way. For any $a \in A$ and $b \in B$, we have
$$f(1_A) = 1_B, \hs \varepsilon_B(f(a)) = \varepsilon_A(a)\hs
\mbox{and}\hs S_B (f(a)) = f(S_A(a)).$$
\\
{\bf 2.4 Proposition} \hs A multiplier Hopf ($^\ast$-)algebra $B$
is $G$-cograded if and only if there exists an injective,
non-degenerate ($^\ast$-)homomorphism $ I : K(G) \rightarrow M(B)$
so that
\begin{itemize}
\item[(i)] $I(K(G)) \subset Z(M(B))$ when $Z(M(B))$ is the center of
$M(B)$,
\item[(ii)] $\Delta(I(f)) = (I \otimes I) \Delta(f)$ for all $f \in
K(G)$.
\end{itemize}
We have $I(\delta_p) = 1_p$, where $\delta_p$ is the complex
valued function on $G$, given by $\delta_p(q) = 0$ if $p\neq q$
and $\delta_p(p)=1$ and where $1_p$ is the unit in $M(B_p)$.
\\ Again, for the proof of
these two results, we refer to [A-D-VD].
\\ \ \\
An immediate
consequence of these two results is the following.
\\
{\bf 2.5 Proposition} \hs Let $B$ be a $G$-cograded multiplier
Hopf $(^\ast$-)algebra in the sense of Definition 2.1. So $B$ has
the form $B = \bigoplus\limits_{p \in G} B_p$.  Then we have
\begin{itemize}
\item[(i)] $\varepsilon (a) = 0$ whenever $a\in B_p$ and $p \neq
e$ (where $e$ is the identity in $G$),
\item[(ii)] $S(B_p) \subseteq M(B_{p^{-1}})$ for all $p$.
\end{itemize}
\ \\
Now that we have discussed the notion of a group-cograded
multiplier Hopf algebra in general, we are ready to study a
special type of actions on these objects.
\\ \ \\
{\it Admissible actions on group-cograded multiplier Hopf
algebras}
\\ \ \\
Here is the main definition.

{\bf 2.6 Definition} \hs Let $B$ be a $G$-cograded multiplier Hopf
($^\ast$-)algebra.  So $B$ has the form $B =
\bigoplus\limits_{p\in G} B_p$.  Let $Aut(B)$ denote the group of
algebra automorphisms on $B$.  By an {\it action} of the group $G$
on $B$, we mean a group homomorphism $\pi : G \rightarrow Aut(B)$.
If $B$ is a multiplier Hopf $^\ast$-algebra we assume that $\pi_p$
is a $^\ast$-automorphism for all $p \in G$.  Further, we require
that for all $p \in G$
\begin{itemize}
\item[(1)] $\pi_p$ respects the comultiplication on $B$ in the
sense that $\Delta(\pi_p(b)) = (\pi_p \otimes \pi_p)(\Delta(b))$
for all $b \in B$.
\end{itemize}
We call this action {\it admissible} if there is an action $\rho$
of $G$ on itself so that
\begin{itemize}
\item[(2)] $\pi_p(B_q) = B_{\rho_p(q)}$,
\item[(3)]  $\pi_{\rho_{p}(q)} = \pi_{pqp^{-1}}$.
\end{itemize}
for all $p,q$ in $G$.
\\ \ \\
The action $\rho$ of $G$ on itself determines an action
$\tilde\rho$ of $G$ on $K(G)$ by the formula $(\tilde\rho_p(f))(q)
= f(\rho_{p^{-1}}(q))$ when $f\in K(G)$ and $p,q\in G$. This is an
action of $G$ on the multiplier Hopf algebra $K(G)$ (in the sense
of the above definition). Condition (2) says that
$I\circ\tilde\rho_p=\pi_p\circ I$ for all $p\in G$ where $I$ is
the canonical imbedding of $K(G)$ in $M(B)$. So, for an action to
be admissible, we first of all need that, on the level of $K(G)$,
it comes from an action of $G$ on itself.
\\
If this action is the adjoint action, that is, if
$\rho_p(q)=pqp^{-1}$ for all $p,q$, then condition (3) is
automatically fulfilled. If this is not the case, then we want
$\pi$ itself to take care of $\rho$ not being the adjoint action.
The other extreme therefore is obtained when $\pi$ is simply the
trivial action. Then we can take for $\rho$ also the trivial
action in order to satisfy (1) while (2) is again automatically
satisfied, now however for a completely different reason.
\\
Condition (2) seems to be quite natural. We will indicate further
why we need condition (3). In any case, as we saw above, we have
the following example.
\\ \ \\
{\bf 2.7 Example} \hs Let $B$ be a $G$-cograded multiplier Hopf
algebra with canonical decomposition $B=\sum_{p\in G}\oplus B_p$.
Let $\pi$ be an action of $G$ on $B$. If $\pi_p(B_q)=B_{pqp^{-1}}$
for all $p,q\in G$, then we have an admissible action.
\\
Also the trivial action is admissible and it is not hard to
construct examples combining these two extreme cases.
\\ \ \\
{\it Deformation of a group-cograded mutliplier Hopf algebra}
\\ \ \\
Let $G$ be a group.  Let $B$ be a $G$-cograded regular multiplier
Hopf algebra and let $\pi$ be an admissible action of $G$ on $B$.
We will construct a new regular multiplier Hopf algebra on $B$ by
deforming the comultiplication, while the algebra structure on $B$
is kept.  The deformation of the comultiplication of $B$, as
defined in the following definition, depends
on the action $\pi$.
\\ \ \\
{\bf 2.8 Definition} \hs Take $B$ and $\pi$ as above.  For $b \in
B$, we define the multiplier $\widetilde{\Delta}(b)$ in $M(B
\otimes B)$ by the following formulas.  Take $b' \in B_q$, then we
define
\begin{eqnarray*}
\widetilde{\Delta}(b) (1 \otimes b') &=& (\pi_{q^{-1}} \otimes
\iota)
(\Delta(b) (1 \otimes b'))\ \\
(1 \otimes b') \widetilde{\Delta}(b) &=& (\pi_{q^{-1}} \otimes
\iota)
 ((1 \otimes b') \Delta(b))
\end{eqnarray*}
in $B\otimes B$.  By the associativity of the product in $B
\otimes B$, we have that
$\widetilde{\Delta}(b)$ is a multiplier in $M(B \otimes B)$.\\
\ \\
{\bf 2.9 Proposition} \hs Take the notations as in Definition 2.8.
For all $b \in B_p$ and $b' \in B_q$, we have
\begin{eqnarray*}
(b' \otimes 1)\widetilde{\Delta}(b) &=& \textstyle\sum b'
\pi_{qp^{-1}}
(b_{(1)}) \otimes b_{(2)} \\
\widetilde{\Delta} (b) (b' \otimes 1) &=& \textstyle\sum
\pi_{qp^{-1}} (b_{(1)}) b' \otimes b_{(2)}
\end{eqnarray*}
in $B\otimes B$.
The map $\widetilde{\Delta} : B \rightarrow M(B \otimes B)$ is a
non-degenerate homomorphism.\\
\ \\
{\bf Proof.} To prove the first formula, take $r\in G$ and $b''
\in B_r$.  Then we have $$((b' \otimes 1) \widetilde{\Delta}(b))
(1 \otimes b'') = (b' \otimes 1) (\widetilde{\Delta}(b) (1 \otimes
b'')) = \textstyle\sum b' \pi_{r^{-1}} (b_{(1)}) \otimes b_{(2)}
b''.$$ As $b' \in B_q$, we must have $\pi_{r^{-1}}(b_{(1)})\in
B_q$ and so $b_{(1)}\in \pi_r(B_q)=B_{\rho_r(q)}$ (using condition
(2) in Definition 2.6). As $b''\in B_r$ we must have $b_{(2)}\in
B_r$. Finally, because $b\in B_p$ we need $p=\rho_r(q)r$. It
follows that $\pi_{pr^{-1}}=\pi_{\rho_r(q)} = \pi_{rqr^{-1}}$
(using condition (3) in Definition 2.6). Therefore,
$\pi_p=\pi_{rq}$ and $\pi_{r^{-1}}=\pi_{qp^{-1}}$. This proves the
first formula.
Similarly, the second one can be proven. \\
Further, $\widetilde{\Delta} : B \rightarrow M(B \otimes B)$ is a
homomorphism because $\Delta : B \rightarrow M(B \otimes B)$ is a
homomorphism and $(\pi_p \otimes \iota_B)$ is a non-degenerate
homomorphism on $B \otimes B$ for all $p \in G$.  For all $p, q
\in G$ we have that $B_p \otimes B_q = \widetilde{\Delta}
(B_{\rho_q(p)q}) (1 \otimes B_q)$ and $B_p \otimes B_q = (1
\otimes B_q) \widetilde{\Delta} (B_{\rho_q(p)q})$.  Therefore,
$\widetilde{\Delta} : B \rightarrow M(B \otimes B)$ is a
non-degenerate homomorphism in the sense that $\widetilde{\Delta}
(B) (B \otimes B) = B \otimes B = (B \otimes B)
\widetilde{\Delta}(B)$. This completes the proof.
 \insp $\blacksquare$
\\ \ \\
{\bf 2.10 Lemma} \hs Take the notations as above.
Then $\widetilde{\Delta}$ is coassociative.\\
\ \\
{\bf Proof.}  As $\widetilde{\Delta}$ can be extended to $M(B)$ in
a natural way, to show that $\widetilde{\Delta}$ is coassociative
we need  $(\widetilde{\Delta} \otimes \iota)
(\widetilde{\Delta}(x)) \stackrel{(\ast)}{=} (\iota \otimes
\widetilde{\Delta}) (\widetilde{\Delta}(x))$ in
$M(B \otimes B \otimes B)$ for all $x \in B$ .\\
Let $1_p$ and $1_q$ denote the units in $M(B_p)$ and $M(B_q)$
respectively. Then the equation $(\ast)$ will be satisfied if for
all $p, q\in G$ we have
\begin{eqnarray*}
((\widetilde{\Delta} \otimes \iota) (\widetilde{\Delta}(x))) (1
\otimes 1_p\otimes 1_q) = ((\iota \otimes \widetilde{\Delta})
(\widetilde{\Delta}(x))) (1 \otimes 1_p\otimes 1_q).
\end{eqnarray*}
For the left hand side of the above equation we set
\begin{eqnarray*}
((\widetilde{\Delta} \otimes \iota) (\widetilde{\Delta}(x))) (1
\otimes 1_p\otimes 1_q) &=& ((\widetilde{\Delta} \otimes \iota)
(\widetilde{\Delta}(x)
(1\otimes 1_q))) (1 \otimes 1_p\otimes 1_q)\\
&=& ((\widetilde{\Delta} \otimes \iota) (\pi_{q^{-1}} \otimes
\iota)
 (\Delta(x))) (1\otimes 1_p\otimes 1_q)\\
&=& ((\pi_{p^{-1}} \otimes \iota \otimes \iota)((\Delta \otimes
\iota)
(\pi_{q^{-1}} \otimes \iota) (\Delta (x))))(1 \otimes 1_p\otimes 1_q)\\
&=& ((\pi_{p^{-1}} \pi_{q^{-1}} \otimes \pi_{q^{-1}} \otimes
\iota) ((\Delta \otimes \iota)(\Delta(x)))) (1\otimes 1_p\otimes
1_q).
\end{eqnarray*}
Observe that $1_p \otimes 1_q = \widetilde{\Delta}
(1_{\rho_q(p)q})(1_p \otimes 1_q)$ because $\pi_{q^{-1}}
(B_{\rho_q(p)}) = B_p$. For the right hand side of the above
equation we set
\begin{eqnarray*}
((\iota \otimes \widetilde{\Delta}) (\widetilde{\Delta}(x))) (1
\otimes 1_p\otimes 1_q) &=& ((\iota \otimes
\widetilde{\Delta})(\widetilde{\Delta}(x))) (1 \otimes
\widetilde{\Delta}(1_{\rho_{q^{(p)}}q})) (1\otimes 1_p\otimes
1_q)\\
&=& ((\iota \otimes \widetilde{\Delta}) (\widetilde{\Delta}(x) (1
\otimes 1_{\rho_q(p)q}))) (1 \otimes 1_p\otimes 1_q)\\
&=& ((\iota \otimes
\widetilde{\Delta})((\pi_{q^{-1}\rho_{q}(p^{-1})} \otimes \iota)
(\Delta(x)))) (1 \otimes 1_p\otimes 1_q).
\end{eqnarray*}
As $\pi$ is an admissible action of $G$ on $B$, we have
$\pi_{q^{-1} \rho_q(p^{-1})} = \pi_{p^{-1} q^{-1}}$.  Therefore,
we obtain that the right hand side equals
\begin{eqnarray*}
\begin{array}{l}\hspace*{-1.0cm}
((\pi_{p^{-1}} \pi_{q^{-1}} \otimes \iota \otimes
\iota) (( \iota \otimes \widetilde{\Delta})(\Delta(x)))) (1
\otimes 1_p\otimes 1_q)\\
\hspace*{1.0 cm}= ((\pi_{p^{-1}} \pi_{q^{-1}} \otimes \iota
\otimes \iota) ((\iota \otimes \pi_{q^{-1}} \otimes \iota)
((\iota \otimes \Delta) (\Delta(x))))) (1 \otimes 1_p\otimes 1_q)\\
\hspace*{1.0 cm}= ((\pi_{p^{-1}} \pi_{q^{-1}} \otimes \pi_{q^{-1}}
\otimes \iota) ((\iota \otimes \Delta)(\Delta(x)))) (1 \otimes
1_p\otimes 1_q).
\end{array}
\end{eqnarray*}
We see that both expressions are the same.  \insp $\blacksquare$\\
\ \\
We now prove that the comultiplication $\widetilde{\Delta}$ makes
$B$ into a regular multiplier Hopf algebra.  We also calculate the
counit and the antipode for this new multiplier Hopf algebra.\\
\ \\
{\bf 2.11 Theorem} \hs Let $B$ be a regular $G$-cograded
multiplier Hopf algebra.  So, as an algebra, $B$ has the form $B =
\bigoplus\limits_{p \in G} B_p$.  Assume that $\pi$ is an
admissible action of $G$ on $B$.  Let $\widetilde{\Delta}$ denote
the comultiplication as defined in Definition 2.8.  Then we have
the following.
\begin{itemize}
\item[(1)] $(B,\widetilde{\Delta})$ is a regular multiplier Hopf
algebra.  The counit $\widetilde{\varepsilon}$ is the original
counit $\varepsilon$. The antipode $\widetilde{S}$ is given by the
formula $\widetilde{S}(b) = \pi_{p^{-1}} (S(b))$ for $b \in B_p$.
\item[(2)] If $B$ is a $G$-cograded multiplier Hopf
$^\ast$-algebra, then $(B, \widetilde{\Delta})$ is again a
multiplier Hopf $^\ast$-algebra.
\item[(3)] If $\varphi$ is a left integral on $B$, then $\varphi$
is also a left integral on $(B,\widetilde{\Delta})$.  However, if
$\psi$ is a right integral, it is in general not right invariant
on $(B,\widetilde{\Delta})$.  It has to be modified.  For $b \in
B_p$, define $\widetilde{\psi} (b) = \psi_B(\pi_{p^{-1}} (b))$.
Then $\widetilde\psi$ is a right integral. In the $^\ast$-case, we
have that a positive left integral on $B$ is again a positive left
integral on $(B, \widetilde{\Delta})$. A positive right integral
$\psi$ on $B$ gives rise to a positive right integral
$\widetilde{\psi}$ on $(B, \widetilde{\Delta})$.
\end{itemize}

{\bf Proof.}
\begin{itemize}
\item[(1)] We will make use of [VD2, Proposition 2.9] to prove
that $(B,\widetilde{\Delta})$ is a regular multiplier Hopf
algebra.  Recall that, as an algebra, $B = \bigoplus\limits_{p\in
G} B_p$.  The comultiplication $\widetilde{\Delta}: B \rightarrow
M(B \otimes B)$ is defined as in Definition 2.8 and we will also
use the formulas of Proposition 2.9.\\
Moreover, we will use that $\pi_p$ is an isomorphism of $B$ for
all $p \in G$ which respects the comultiplication in the sense
that $\Delta(\pi_p(b)) = (\pi_p \otimes \pi_p) (\Delta(b))$ for
all $b \in B$.\\
\ \\
We first consider the counit. We will show that the original
counit $\varepsilon$ on $B$ is also the counit for $\widetilde B$.
Take $b \in B_p$ and $b' \in B_q$.  Then we have
\begin{eqnarray*}
\begin{array}{l}
(\varepsilon \otimes \iota) (\widetilde{\Delta} (b) (1 \otimes
b')) = \sum \varepsilon(\pi_{q^{-1}} (b_{(1)})) b_{(2)}b' =
\sum\varepsilon (b_{(1)})b_{(2)} b' = bb'
\\
(\iota \otimes \varepsilon)((b'\otimes 1) \widetilde{\Delta}(b)) =
(\iota \otimes \varepsilon) ((b' \otimes 1) \widetilde{\Delta}(b)
(1 \otimes 1_e)) =
\\
\hspace*{1cm}(\iota \otimes \varepsilon) ((b' \otimes 1) (\pi_e
\otimes \iota) (\Delta(b) (1\otimes 1_e))) = (\iota \otimes
\varepsilon) ((b' \otimes 1) \Delta(b)) = b'b.
\end{array}
\end{eqnarray*}
Recall that $e$ denotes the identity in $G$. \\ Next, we prove the
existence of the antipode. Define $\widetilde{S}$ by the formula
$\widetilde{S}(b) = \pi_{p^{-1}} (S(b))$ for all $b \in B_p$. Let
$m$ denote the multiplication in the algebra $B$. Take $b \in B_p$
and $b' \in B_q$. Then we have
\begin{eqnarray*}
m((\widetilde{S} \otimes \iota)(\widetilde{\Delta}(b) (1 \otimes
b'))) = m(\widetilde{S} \otimes \iota) (\textstyle\sum
\pi_{q^{-1}} (b_{(1)}) \otimes b_{(2)} b').
\end{eqnarray*}
As $b_{(2)} \in B_q$, we calculate that $\pi_{q^{-1}} (b_{(1)})
\in B_{\rho_{q^{-1}} (pq^{-1})}$.  Notice that $\pi_{\rho_{q^{-1}}
(qp^{-1})} = \pi_{q^{-1}} \pi_{q\rho_{q^{-1}}(qp^{-1})} =
\pi_{p^{-1}q}$.\\
Therefore, we get
\begin{eqnarray*}
m((\widetilde{S} \otimes \iota) (\widetilde{\Delta}(b) (1 \otimes
b'))) = \textstyle\sum \pi_{p^{-1}} (S(b_{(1)})) b_{(2)} b'.
\end{eqnarray*}
If $\rho_{p^{-1}} (qp^{-1}) \neq q$ (and hence $p \neq e$), then
the last expression equals zero. Remark that also
$\widetilde{\varepsilon} (b) b' = 0$. If $\rho_{p^{-1}} (qp^{-1})
= q$, then $\pi_q = \pi_{\rho_{p^{-1}}(qp^{-1})} = \pi_{p^{-1}q}$.
Therefore we have that in this case $\pi_{p^{-1}} = \pi_e$ and the
expression in the right hand side becomes $\sum S(b_{(1)}) b_{(2)}
b' = \varepsilon (b) b' = \widetilde{\varepsilon}(b)
b'$.
\\We now
prove the second equation for $\widetilde{S}$. Take $b \in B_p$
and $b' \in B_q$, then we have
\begin{eqnarray*}
m((\iota \otimes \widetilde{S}) ((b' \otimes
1)\widetilde{\Delta}(b))) = m((\iota \otimes \widetilde{S})
(\textstyle\sum b' \pi_{qp^{-1}} (b_{(1)}) \otimes b_{(2)})).
\end{eqnarray*}
In this formula we calculate that $b_{(2)} \in
B_{\rho_{pq^{-1}}(q^{-1})p}$. Notice that
$\pi_{p^{-1}\rho_{pq^{-1}}(q)} = \pi_{qp^{-1}}$.\\
Therefore, we obtain
\begin{eqnarray*}
m((\iota \otimes \widetilde{S}) ((b' \otimes
1)\widetilde{\Delta}(b))) = \textstyle\sum b' \pi_{qp^{-1}}
(b_{(1)} S(b_{(2)})) = \widetilde{\varepsilon} (b) b'.
\end{eqnarray*}
By using [VD2, Proposition 2.9], we conclude that $(B,
\widetilde{\Delta})$ is a regular multiplier Hopf algebra.
\item[(2)] Now assume that $B$ is a $G$-cograded multiplier Hopf
$^\ast$-algebra and for all $p \in G$, $\pi_p$ is furthermore a
$^\ast$-isomorphism on $B$. We prove that $(B,\widetilde{\Delta})$
is also a multiplier Hopf $^\ast$-algebra.  Therefore, we have to
show that $\widetilde{\Delta}$ is a $^\ast$-homomorphism.  Take $b
\in B_p$ and $b' \in B_q$, then we have
\begin{eqnarray*}
\widetilde{\Delta}(b^\ast) (1 \otimes b') &=& (\pi_{q^{-1}}
\otimes \iota) (\Delta(b^\ast) (1 \otimes b')) = (\pi_{q^{-1}}
\otimes \iota) (\Delta(b)^\ast (1 \otimes b'))\\
&=& (\pi_{q^{-1}} \otimes \iota) (((1 \otimes b'{^\ast}) \Delta
(b))^\ast) = ((\pi_{q^{-1}} \otimes \iota) ((1 \otimes b'{^\ast})
\Delta (b)))^\ast\\
 &=& ((1 \otimes b'{^\ast})
\widetilde{\Delta}(b))^\ast = \widetilde{\Delta}(b)^\ast (1
\otimes b').
\end{eqnarray*}
\item[(3)] Let $\varphi$ be a left integral on $B$, as reviewed in
Section 1. We prove that $\varphi$ is also a left integral on
$(B,\widetilde{\Delta})$.  Take $b \in B_p$ and $b' \in B_q$.
Then, we have $((\iota \otimes \varphi) \widetilde{\Delta}(b))b' =
(\iota \otimes \varphi) (\widetilde{\Delta}(b) (b' \otimes 1)) =
\sum \pi_{qp^{-1}} (b_{(1)}) b'\varphi (b_{(2)}) = \varphi(b)b'$.
Therefore, $(\iota \otimes \varphi)\widetilde{\Delta}(b) =
\varphi(b) 1$ in $M(B)$ for all $b \in B$.\\
Let $\psi$ be a right integral on $B$.  Define $\widetilde{\psi}$
on $B$ by the formula $\widetilde{\psi}(b) = \psi_B (\pi_{p^{-1}}
(b))$ when $b \in B_p$.  We will now show that $\widetilde{\psi}$
is a right integral on $(B,\widetilde{\Delta})$. Take $b \in B_p$
and $b' \in B_q$, then we have
\begin{eqnarray*}
((\widetilde{\psi} \otimes \iota) \widetilde{\Delta}(b)) b' =
(\widetilde{\psi} \otimes \iota)(\widetilde{\Delta}(b) (1 \otimes
b')) = (\widetilde{\psi} \otimes \iota)
(\textstyle\sum\pi_{q^{-1}} (b_{(1)}) \otimes b_{(2)} b').
\end{eqnarray*}
As $\pi_{q^{-1}} (b_{(1)}) \in B_{\rho_{q^{-1}}(pq^{-1})}$, the
last expression is given as
\begin{eqnarray*}
\textstyle\sum \psi (\pi_{\rho_{q^{-1}}(qp^{-1})} \pi_{q^{-1}}
(b_{(1)}))
b_{(2)} b' &=& \textstyle\sum \psi (\pi_{p^{-1}} (b_{(1)})) b_{(2)} b'\\
&=& \pi_p (\textstyle\sum\psi (\pi_{p^{-1}}(b_{(1)})) \pi_{p^{-1}}
(b_{(2)})
\pi_{p^{-1}}(b'))\\
&=& \pi_p (\psi (\pi_{p^{-1}} (b)) \pi_{p^{-1}}(b')) =
\widetilde{\psi}(b)b'.
\end{eqnarray*}
If $B$ is a $G$-cograded multiplier Hopf $^\ast$-algebra as in
(2), then a positive left integral $\varphi$ on $B$ stays a
positive left integral on $(B,\widetilde{\Delta})$ because the
$^\ast$-algebra structure of $B$ and $(B, \widetilde{\Delta})$ is
the same. If $\psi$ is a positive right integral on $B$, then it
is easily shown that $\widetilde{\psi}$ is a positive right
integral on $(B, \widetilde{\Delta})$.\insp $\blacksquare$\\
\end{itemize}

{\bf 2.12 Notation} \hs We will use $\widetilde{B}$ for the
deformed multiplier Hopf algebra $(B,\widetilde{\Delta})$.  Of
course, if $\pi$ is the trivial action of $G$ on $B$, we have that
$(B,\widetilde{\Delta})$ equals $B$.
\\ \ \\
In Proposition 3.12 we will refine the structure of
$\widetilde{B}$ in the case where $\pi$ is a admissible action
such that $\pi_p(B_q) = B_{pqp^{-1}}$ for all $p, q \in G$.  This
case is like the {\it mirror} construction for a Hopf
group-coalgebra as introduced in [T, Section 11].

\section{Pairing and Drinfel'd double construction with
\boldmath{$G$}-cograded multiplier Hopf ($^\ast$-)algebras} In
this section, we will apply the 'twisted tensor product'
construction of multiplier Hopf ($^\ast$-)algebras as we have
explained in Section 1.
\\
Let $G$ be a group and let $B$ be a regular $G$-cograded
multiplier Hopf algebra in the sense of Definition 2.1.  We
suppose that $\pi$ is an admissible action of $G$ on $B$ in the
sense of Definition 2.6. In Proposition 3.1 below we study a
multiplier Hopf algebra pairing $\langle A, B \rangle$ when $B$ is
$G$-cograded. Further in this section, we define Drinfel'd double
constructions $D^\pi = A^{cop} \bowtie \widetilde{B}$ where the
product as well as the coproduct are depending on the action
$\pi$.\\
\ \\
{\bf 3.1 Proposition}\hs Let $\langle A,B\rangle$ be a pairing of
two regular multiplier Hopf algebras.  Suppose that $B$ is
$G$-cograded.  Then there exist subspaces $\{A_p\}_{p\in G}$ of
$A$ such that
\begin{itemize}
\item[(1)] $A = \bigoplus\limits_{p\in G} A_p$ and
$A_pA_q\subseteq A_{pq}$,
\item[(2)] $\langle A_p, B_q\rangle = 0$ whenever $p\neq q$,
\item[(3)] $\langle \Delta(A_p), B_q \otimes B_r\rangle = 0$
if $q \neq p$ or $r \neq p$,
\end{itemize}
where $p,q,r\in G$. \\
{\bf Proof.} \hs As $B$ is a $G$-cograded
multiplier Hopf algebra, $B$ has the form $B =
\bigoplus\limits_{p\in G} B_p$. Recall from Section 1, that there
are four module algebra structures associated to the pairing
$\langle A, B \rangle$, denoted as $A \btr B$, $B \btr A$, $A \btl
B$, $B \btl A$. We have seen that these actions are unital and
therefore extend to the multiplier algebras. One can show that  $A
\btl B_p= B_p\btr A = A \btl 1_p= 1_p\btr A$ for all $p\in G$
where $1_p$ denotes the unit in $M(B_p)$. Then we define a
subspace $A_p$ in $A$ by $A_p=1_p\btr A$ for any $p\in G$. We now
prove that these subspaces satisfy the 3 requirements.
\begin{itemize}
\item[(1)] Take $a \in A$. Then there exists an element $b \in B$
so that $a = b\btr a$.  As $B = \bigoplus\limits_{p \in G} B_p$,
it easily follows that $A = \bigoplus\limits_{p\in G} A_p$.
Furthermore $A_pA_q = (1_p\btr A_p) (1_q\btr A_q) = 1_{pq}\btr
(A_pA_q)$. Therefore, $A_pA_q\subseteq A_{pq}$.
\item[(2)] This second property follows from the definition of
$A_p$ and the multiplication structure of $B$.
\item[(3)] Take $a \in A_p$, $b \in B_q$ and $b' \in B_r$.  Then
we have that $\langle \Delta(a), b \otimes b'\rangle = \langle a,
bb'\rangle$.  Now the result follows from the algebra structure on
$B$ and the fact that $A_p$ is also given by $A_p=  B_p\btr A$.
\hfill $\blacksquare$
\end{itemize}
\ \\
{\bf 3.2 Remark}\hs Take a pairing $\langle A, B\rangle$ as in
Proposition 3.1.
\begin{itemize}
\item[(1)] It easily follows from the definition that $b\btr a=0$
if $a\in A_q$, $b\in B_p$ and $p\neq q$. Indeed, when $a\in A_q$
and $b\in B$, we have $b\btr a=(b1_q) \btr a$ and $b1_q=0$ if
$b\in B_p$ and $p\neq q$. Furthermore, if $a\in A_q$ then $b\btr
a\in A_q$ for all $b\in B$. Similar results hold for the module
$A\btl B$.
\item[(2)] The antipode $S$ of $A$ maps $A_p$ to $A_{p^{-1}}$.
\end{itemize}

Again take a pairing $\langle A, B  \rangle$ as in Proposition 3.1
and now let $\pi$ be an admissible action of $G$ on $B$.  We will
construct a twisted tensor product multiplier Hopf algebra, as
reviewed in Section 1.  The twist map, defining the non-trivial
product structure on $A \otimes B$, will depend on the pairing as
well as on the action $\pi$.  The comultiplication, which is
compatible with this product on $A \otimes B$, will also depend on
the action $\pi$.
\\ \ \\
We first prove the following lemma.\\
{\bf  3.3 Lemma}\hs Take the notations and the assumptions as
above.   Define linear maps $R_1$ and $R_2$ on $A \otimes B$ by
the formulas
\begin{eqnarray*}
R_1 (a \otimes b) &=& \textstyle\sum (\pi_{qp^{-1}} (b_{(1)}) \btr a) \otimes b_{(2)}\\
R_2 (a \otimes b) &=& \textstyle\sum (a \btl b_{(2)}) \otimes b_{(1)}
\end{eqnarray*}
when $a \in A_q$ and $b \in B_p$. Then $R_1$ and $R_2$ are
bijections on $A \otimes B$ and the inverses are given by
\begin{eqnarray*}
R_1^{-1} (a \otimes b) &=& \textstyle\sum (\pi_{p^{-1}} (S^{-1} (b_{(1)})) \btr a) \otimes b_{(2)}\\
R_2^{-1} (a \otimes b) &=& \textstyle\sum (a \btl S^{-1} (b_{(2)})) \otimes b_{(1)}
\end{eqnarray*}
when $b \in B_p$.\\
{\bf Proof.} \hs We remark that in all the formulas above, the
decompositions are well-covered because the modules $B \btr A$ and
$A \btl B$ are unital.
\\
The proof for the map $R_2$ is easy and the result is known. Here,
we do not really need these restrictions on $a$ and $b$.
\\
We give the proof for the map $R_1$. Take $a \in A_q$ and $b \in
B_p$.  First remark that, when looking closer at the definition of
$R_1^{-1}$, we see that $\pi_{p^{-1}}(S^{-1}(b_{(1)}))$ is forced
to lie in $B_q$ as it acts on the element $a$ in $A_q$; see Remark
3.2 (1). Then $b_{(1)}$ must be in $B_{\rho_p(q^{-1})}$. Because
we assume that $b\in B_p$ it follows that that $b_{(2)}$ is forced
to lie in $B_{\rho_p(q)p}$.  This is used when we apply the map
$R_1$ in the following calculation. We have
\begin{eqnarray*}
(R_1 \circ R_1^{-1})(a\otimes b)
&=& R_1(\sum(\pi_{p^{-1}} (S^{-1} (b_{(1)})) \btr a) \otimes b_{(2)})\\
&=& \sum ((\pi_{qp^{-1}\rho_{p}(q^{-1})} (b_{(2)}) \pi_{p^{-1}} (S^{-1} (b_{(1)}))  \btr a) \otimes b_{(3)}\\
&=& \sum (\pi_{p^{-1}} (b_{(2)} S^{-1}(b_{(1)})) \btr a) \otimes
b_{(3)} = a \otimes b.
\end{eqnarray*}
Remark that we have used that
$\pi_{qp^{-1}\rho_p(q^{-1})}=\pi_{qp^{-1}pq^{-1}p^{-1}}=\pi_{p^{-1}}$
in the above calculation. A similar argument will give
\begin{eqnarray*}
(R_1^{-1} \circ R_1) (a \otimes b)
&=& R_1^{-1} (\sum(\pi_{qp^{-1}} (b_{(1)}) \btr a) \otimes b_{(2)})\\
&=& \sum ((\pi_{p^{-1} \rho_{pq^{-1}}(q)} (S^{-1} (b_{(2)})) \pi_{qp^{-1}} (b_{(1)})) \btr a) \otimes b_{(3)}\\
&=& \sum (\pi_{qp^{-1}} (S^{-1} (b_{(2)})b_{(1)}) \btr a) \otimes
b_{(3)} = a \otimes b. \insp \blacksquare
\end{eqnarray*}
\ \\
We now define the twist map $R : B \otimes A \rightarrow A \otimes B$.
\\
{\bf  3.4 Definition}\hs Take a pairing $\langle A, B\rangle$ as
in Proposition 3.1.  Let $\pi$ be an admissible action of $G$ on
$B$.  We define the twist map $R : B \otimes A \rightarrow A
\otimes B$ by the composition $R = R_1 \circ R_2^{-1} \circ
\sigma$, where $\sigma$ is the flip map from $B\otimes A$ to
$A\otimes B$. So, for all $a \in A_q$ and $b \in B_p$ we have,
using arguments as in the proof of the previous lemma,
$$R(b \otimes a) = \textstyle\sum (\pi_{p^{-1}} (b_{(1)}) \btr a \btl S^{-1} (b_{(3)})) \otimes b_{(2)}.$$
As a composition of bijections, $R$ is a bijection.
\\ \ \\
We see from this formula, that it will be true for all $a\in A$
when $b\in B_p$. When $a\in A_q$ and $b\in B_p$ then $b_{(2)}$ is
forced in $B_{\rho_p(q^{-1})pq}$ in the above formula. This will
be used, in particular, in the proof of the following lemma where
we obtain that $R$ behaves well with respect to the
multiplications of $A$ and $B$.
\\
{\bf 3.5 Lemma}\hs
Take the notations and assumptions as before.
The twist map $R$ satisfies the following equations
\begin{eqnarray*}
\begin{array}{llll}
(1) &R  (m_B \otimes \iota_A)
= (\iota_A \otimes m_B)(R \otimes \iota_B) (\iota_B \otimes R) &\mbox{\quad on \quad} &B\otimes B \otimes A,\\
(2) &R  (\iota_B \otimes m_A) = (m_A \otimes \iota_B)(\iota_A
\otimes R) (R \otimes \iota_A) &\mbox{\quad on \quad} &B \otimes A
\otimes A,
\end{array}
\end{eqnarray*}
where, as before, $m_A$ and $m_B$ are the multiplications in $A$
and $B$ and where $\iota_A$ and $\iota_B$ are the identity maps on
$A$ and $B$ respectively.\\
{\bf Proof.} (1) Take $b \in B_p$, $b' \in B_q$ and $a \in A_r$.
Let the right hand side of the first equation act on $b \otimes b'
\otimes a$. We get
\begin{eqnarray*}
\begin{array}{l}
((\iota_A \otimes m_B) (R \otimes \iota_B)(\iota_B \otimes R)) (b \otimes b' \otimes a)\\
\hspace*{2 cm}= \sum(\pi_{p^{-1}} (b_{(1)}) \pi_{q^{-1}}
(b'_{(1)}) \btr a \btl S^{-1} (b'_{(3)}) S^{-1} (b_{(3)})) \otimes
b_{(2)} b'_{(2)}.
\end{array}
\end{eqnarray*}
By the remark preceding this lemma, we find that $b_{(2)} \in
B_{\rho_p (r^{-1}) pr}$ and $b'_{(2)} \in B_{\rho_q(r^{-1}) qr}$.
Therefore, $\rho_p(r^{-1}) pr = \rho_q(r^{-1}) qr$ and we obtain
that $\pi_p = \pi_q$.  The above equation can now be written as
\begin{eqnarray*}
\begin{array}{l}
((\iota_A \otimes m_B)(R \otimes \iota_B)(\iota_B \otimes R)) (b\otimes b' \otimes a)\\
\hspace*{2 cm}= \sum (\pi_{p^{-1}} (b_{(1)} b'_{(1)}) \btr a \btl
S^{-1} (b_{(3)} b'_{(3)})) \otimes b_{(2)} b'_{(2)}.
\end{array}
\end{eqnarray*}
\begin{itemize}
\item[$\bullet$]  If $p \neq q$, this expression equals zero
because $bb'=0$. Clearly the operator $R (m_B \otimes \iota_A)$
applied on $(b \otimes b' \otimes a)$ also equals zero in this
case.
\item[$\bullet$] If $p = q$, then $b, b' \in B_p$ and also $bb' \in B_p$.
The above equation now becomes
\begin{eqnarray*}
((\iota_A \otimes m_B)(R \otimes \iota_B) (\iota_B \otimes R)) (b
\otimes b' \otimes a) = (R  (m_B \otimes \iota_A)) (b \otimes b'
\otimes a).
\end{eqnarray*}
\end{itemize}
This completes the proof of the first statement.\\
(2) To prove the second statement, take $b \in B_p$, $a \in A_q$
and $a' \in A_{r}$. Let the right hand side of the second equation
act on $b \otimes a \otimes a'$.  We get
\begin{eqnarray*}
\begin{array}{l}
((m_A \otimes \iota_B)(\iota_A \otimes R) (R \otimes \iota_A))
(b \otimes a \otimes a')\\
\hspace*{2cm}= (m_A \otimes \iota_B) (\iota_A \otimes
R)(\sum(\pi_{p^{-1}} (b_{(1)}) \btr a \btl S^{-1} (b_{(3)}))
\otimes b_{(2)} \otimes a').
\end{array}
\end{eqnarray*}
Recall that $b_{(2)} \in B_{\rho_p(q^{-1}) pq}$. We also have
$\pi_{q^{-1} p^{-1} \rho_p(q)} = \pi_{p^{-1}}$. Therefore, the
above expression can be written as
\begin{eqnarray*}
\begin{array}{l}
(m_A \otimes \iota_B)(\sum(\pi_{p^{-1}}(b_{(1)}) \btr a \btl
S^{-1}(b_{(5)})) \otimes (\pi_{p^{-1}} (b_{(2)}) \btr a' \btl
S^{-1} (b_{(4)})) \otimes b_{(3)})\\
 \hspace*{2.5cm}=
\sum(\pi_{p^{-1}}(b_{(1)}) \btr a \btl S^{-1} (b_{(5)})
(\pi_{p^{-1}} (b_{(2)}) \btr a' \btl S^{-1}(b_{(4)})) \otimes
b_{(3)}.
\end{array}
\end{eqnarray*}
As $B \btr A$ and $A \btl B$ are module algebras, this expression equals
\begin{eqnarray*}
\begin{array}{l}
\sum(\pi_{p^{-1}} (b_{(1)}) \btr (aa') \btl S^{-1} (b_{(3)}))
\otimes b_{(2)}= R(b \otimes aa') = (R  (\iota_B \otimes m_A)) (b
\otimes a \otimes a'). \insp \blacksquare
\end{array}
\end{eqnarray*}
\ \\
As reviewed in Section 1 (Twisted tensor product construction of
multiplier Hopf \\($^*$-)algebras), the  map $R$ defines a
non-trivial product on $A \otimes B$ which is non-degenerate. The
algebra defined in this way is denoted as $A \bowtie B$. Recall
that the product in $A \bowtie B$ is given by the formula
$$(a \bowtie b) (a' \bowtie b') = (m_A \otimes m_B) (\iota_A \otimes R
\otimes \iota_B) (a \otimes b \otimes a' \otimes b')$$ for all $a,
a' \in A$ and $b, b' \in B$. In Section 1, Remarks 1.11, other
expressions are given for the right hand side.
\\ \ \\
We now consider these algebras with their comultiplications. Let
$A^{cop} = (A,\Delta^{cop})$ and $\widetilde{B} = (B,
\widetilde{\Delta})$. In Definition 1.12, we saw that
$\Delta^{cop}(a) \widetilde{\Delta}(b)$ is a multiplier in $M((A
\bowtie B) \otimes (A \bowtie B))$ for all $a \in A$ and $b \in
B$ and that we get a comultiplication. So, the following definition
is possible here.\\
\ \\
{\bf 3.6 Definition}\hs Take the notations and assumptions as
above.  For $a \in A$ and $b \in B$, we define the multiplier
$\overline{\Delta} (a \bowtie b)$ in $M((A \bowtie B) \otimes (A
\bowtie B))$ by the formula
$$\overline{\Delta} (a \bowtie b) = \Delta^{cop} (a) \widetilde{\Delta}(b).$$
\ \\
As we reviewed in Theorem 1.13, for $\overline{\Delta}$ to be a
homomorphism on $A \bowtie B$, we need to prove the following
compatibility relation between
$R$ and $\overline{\Delta}$.\\

{\bf 3.7 Proposition} \hs Take the notations and assumptions as
above.  Then, we have
$$\overline{\Delta} (R(b \otimes a)) = \widetilde{\Delta}(b) \Delta^{cop}(a)$$
in $M((A \bowtie B) \otimes (A \bowtie B))$ for all $a \in A$ and $b \in B$.\\
\ \\
{\bf Proof.} Take $a \in A_{q}$, $a' \in A_{q'}$, and $b \in
B_{p}$, $b' \in B_{p'}$.  Then the product in the twisted tensor
product algebra $A \bowtie B$ is given by the formula
$$(a \bowtie
b) (a' \bowtie b') = \textstyle\sum \langle a'_{(1)}, S^{-1}
(b_{(3)}\rangle \langle a'_{(3)}, \pi_{p^{-1}} (b_{(1)})\rangle
aa'_{(2)} \bowtie b_{(2)} b'.$$ Observe that in the right hand
side, all the decompositions  are covered. Recall also that
$\overline{\Delta} (a \bowtie b) = \Delta^{cop} (a)
\widetilde{\Delta}(b).$
\\
Take $a \in A_q$, $a' \in A_{q'}$, $a'' \in A_{q''}$, $b \in B_p$,
$b'\in B_{p'}$, $b'' \in B_r$ and $y \in B_s$. Then we calculate
in $(A \bowtie B) \otimes (A \bowtie B)$ that
\begin{eqnarray*}
\begin{array}{l}
((a' \bowtie 1) \otimes (a'' \bowtie y)) (\widetilde{\Delta}(b)
\Delta^{cop}(a)) ((1 \bowtie b') \otimes (1 \bowtie b''))\\
\hspace*{2 cm}=\sum((a' \bowtie \pi_{s^{-1}} (b_{(1)}) \otimes
(a'' \bowtie
yb_{(2)})) ((a_{(2)} \bowtie b') \otimes (a_{(1)} \bowtie b''))\\
\hspace*{2cm}= \sum ((a' \bowtie \pi_{s^{-1}} (b_{(1)})) (a_{(2)}
\bowtie b')) \otimes ((a'' \bowtie yb_{(2)})(a_{(1)} \bowtie b''))
\end{array}
\end{eqnarray*}
Now, observe that $\pi_{s^{-1}} (b_{(1)}) \in
B_{\rho_{s^{-1}}(ps^{-1})}$ and $\pi_{\rho_{s^{-1}}(sp^{-1})} =
\pi_{p^{-1}s}$. Then we apply the commutation rules to commute the
the elements $\pi_{s^{-1}}(b_{(1)})$ and $a_{(2)}$ in the first
factor of the tensor product and the elements $b_{(2)}$ and
$a_{(1)}$ in the second factor. Using the property of the
antipode, we finally obtain that the above expression equals
$$\sum\langle a_{(4)}, \pi_{p^{-1}} (b_{(1)})\rangle \langle S^{-1}
(a_{(1)}), b_{(4)}\rangle(a' a_{(3)} \bowtie \pi_{s^{-1}}
(b_{(2)}) b') \otimes ((a'' \bowtie y) (a_{(2)} \bowtie b_{(3)}
b'')).$$ In the second factor of this tensor product, we deal with
a product in $A \bowtie B$. This product equals zero if $r \neq
\rho_s (q^{-1}) qs$.  If $r = \rho_s(q^{-1}) sq$, then
$\pi_{r^{-1}} = \pi_{s^{-1}}$.\\ On the other hand, we also
calculate in $(A \bowtie B) \otimes (A \bowtie B)$ that
\begin{eqnarray*}
\begin{array}{l}
((a' \bowtie 1) \otimes (a'' \bowtie y)) \overline{\Delta} (R(b
\otimes a)) ((1 \bowtie b') \otimes (1 \bowtie b''))\\
\hspace*{.5cm}= ((a' \bowtie 1) \otimes (a'' \bowtie y))
\overline{\Delta} (\sum(\pi_{p^{-1}} (b_{(1)}) \btr a \btl
S^{-1}(b_{(3)})) \bowtie
b_{(2)}) ((1 \bowtie b') \otimes (1 \bowtie b''))\\
\hspace*{.5cm}= \sum((a'\bowtie 1) \otimes (a'' \bowtie y))
\Delta_A^{cop} (\pi_{p^{-1}} (b_{(1)}) \btr a \btl S^{-1}
(b_{(3)})) \widetilde{\Delta}
(b_{(2)}) ((1 \bowtie b') \otimes (1 \bowtie b''))\\
\hspace*{.5cm}= \sum(a'(\pi_{p^{-1}} (b_{(1)}) \btr a_{(2)})
\bowtie \pi_{r^{-1}} (b_{(2)})b') \otimes ((a'' \bowtie
y)((a_{(1)} \btl
S^{-1} (b_{(4)})) \bowtie b_{(3)} b'')\\
\hspace*{.5cm}= \sum \langle a_{(1)}, S^{-1} (b_{(4)}) \rangle
\langle a_{(4)}, \pi_{p^{-1}} (b_{(1)})\rangle(a' a_{(3)} \bowtie
\pi_{r^{-1}} (b_{(2)})b') \otimes ((a'' \bowtie y)(a_{(2)} \bowtie
b_{(3)} b'')).
\end{array}
\end{eqnarray*}
As before, we have a product in $A \bowtie B$ in the second factor
of this tensor product and if $r \neq \rho_s(q^{-1}) sq$, this
last expression equals zero. If
$r = \rho_s(q^{-1}) sq$, we have $\pi_{r^{-1}} = \pi_{s^{-1}}$.\\
As both calculations give the same result in $(A \bowtie B)
\otimes (A  \bowtie B)$, we conclude that $\widetilde{\Delta}(b)
\Delta^{cop} (a) = \overline{\Delta} (R(b \otimes a))$ in $M((A
\bowtie B) \otimes (A \bowtie B))$. \hfill
$\blacksquare$\\
\ \\
We now formulate the main result of this section.\\
{\bf  3.8 Theorem}\hs Let $\langle A, B\rangle$ be a pair of
multiplier Hopf algebras and assume that $B$ is a (regular)
$G$-cograded multiplier Hopf algebra. Let $\pi$ be an admissible
action of $G$ on $B$.
\begin{itemize}
\item[(1)] The space $D^{\pi} = A^{cop} \bowtie \widetilde{B}$ becomes a
(regular) multiplier Hopf algebra, called the Drinfel'd double,
with the multiplication, the comultiplication, the counit and the
antipode, depending on the pairing as well as on the action $\pi$,
defined in the following way:
\begin{itemize}
\item[$\bullet$] $(a \bowtie b) (a' \bowtie b') = (m_A \otimes
m_B) (\iota_A\otimes R \otimes \iota_B) (a \otimes b \otimes a'
\otimes b')$ where $R(b \otimes a') = \sum (\pi_{p^{-1}} (b_{(1)})
\btr a' \btl S^{-1} (b_{(3)})) \otimes b_{(2)}$ for all $a' \in A$
and $b \in B_p$,
\item[$\bullet$] $\overline{\Delta} (a \bowtie b) = \Delta^{cop}
(a) \widetilde{\Delta}(b)$ for all $a \in A$ and $b \in B$ where
$\Delta^{cop}(a)$ and $\widetilde{\Delta}(b)$ are considered as
multipliers in $M(D^\pi \otimes D^\pi)$,
\item[$\bullet$] $\overline{\varepsilon} (a\bowtie b) =
\varepsilon(a) \varepsilon (b)$ for all $a \in A$ and $b \in B$,
\item[$\bullet$] $\overline{S} (a \bowtie b) = R(\pi_{p^{-1}}
(S(b)) \otimes S^{-1}(a))$ for all $a \in A$ and $b \in B_p$.
\end{itemize}
\item[(2)] If moreover $A$ is a multiplier Hopf $^\ast$-algebra, $B$ a
$G$-cograded multiplier Hopf $^\ast$-algebra and $\langle A, B
\rangle$ a $^\ast$-pairing, then $D^\pi$ is again a multiplier
Hopf $^\ast$-algebra with the $^\ast$-operation  given by $(a
\bowtie b)^\ast = R(b^\ast \otimes a^\ast)\ \mbox{for all $a \in
A$ and $b \in B$.}$
\end{itemize}

{\bf Proof.}
\begin{itemize}
\item[(1)] Let $R$ be the twist map, defined in Definition 3.4.
Recall that for $a \in A$ and $b \in B_p$, $R$ is given by the
formula $$R(b \otimes a) = \Sigma(\pi_{p^{-1}} (b_{(1)}) \btr a
\btl S^{-1} (b_{(3)})) \otimes b_{(2)}.$$ This twist map $R : B
\otimes A \rightarrow A \otimes B$ is bijective and satisfies the
compatibility conditions with the multiplications of the algebras
$A$ and $B$, see Lemma 3.5.  As reviewed in Section 1, we consider
the twisted tensor product algebra $A \bowtie B$ associated with
this twist map $R$.  For all $a \in A$ and $b \in B$, we consider
the multiplier $\overline{\Delta} (a \bowtie b) = \Delta^{cop} (a)
\widetilde{\Delta}(b)$ in $M((A \bowtie B) \otimes (A \bowtie
B))$.  Then $\overline{\Delta}$ satisfies the compatibility
condition with $R$ as proven in Proposition 3.7. Following Theorem
1.13 we can consider the twisted tensor product multiplier Hopf
algebra associated to $A^{cop}$, $\widetilde{B}$ and the twist map
$R$. We denote this multiplier Hopf algebra as $D^\pi = A^{cop}
\bowtie \widetilde{B}$. The counit and the antipode on $D^\pi$ are
uniquely determined in this setting and the formulas are given in
the formulation of Theorem 1.13. Remember that the antipode in
$A^{cop}$ is $S^{-1}$ while the antipode in $\widetilde B$ is
given by $\widetilde S(b)=\pi_{p^{-1}}(S(b))$ when $b\in B_p$, see
Theorem 2.11.
\item[(2)] From the conditions on the multiplier Hopf
$^\ast$-algebra $B$, we have that the multiplier Hopf algebra
$\widetilde{B} = (B, \widetilde{\Delta})$ is again a multiplier
Hopf $^\ast$-algebra, see Theorem 2.11 (2).  Following Theorem
1.13 the twisted tensor product $D^\pi = A^{cop} \bowtie
\widetilde{B}$ is a multiplier Hopf $^\ast$-algebra via the
formula $(a \bowtie b)^\ast = R(b^\ast \otimes a^\ast)$ if this
operation defines an involution on $D^\pi$.  To show that this is
the case, take $a \in A_q$ and $b \in B_p$.  Then we have that
\begin{eqnarray*}
\begin{array}{l}
((a \bowtie b)^\ast)^\ast = (R(b^\ast \otimes a^\ast))^\ast =
(\sum(\pi_{p^{-1}} (b^\ast_{(1)}) \btr a^\ast \btl S^{-1}
(b^\ast_{(3)})) \bowtie b^\ast_{(2)})^\ast\\
\hspace*{5cm}= \sum R(b_{(2)} \otimes (\pi_{p^{-1}} (b_{(1)}^\ast)
\btr a^\ast \btl S^{-1} (b^\ast_{(3)}))^\ast).
\end{array}
\end{eqnarray*}
Because $\langle A,B\rangle$ is a $^\ast$-pairing of multiplier
Hopf $^\ast$-algebras, we have that\\ $(b \btr a \btl b')^\ast =
S^{-1} (b^\ast) \btr a^\ast \btl S^{-1}(b'{^\ast})$ for all $a \in
A$ and $b, b'\in B$.  Therefore, we have $((a \bowtie
b)^\ast)^\ast = \sum R(b_{(2)} \otimes (\pi_{p^{-1}} (S^{-1}
(b_{(1)})) \btr a \btl b_{(3)}))$.  Because $a \in A_q$, we must
have that $b_{(2)} \in B_{\rho_p(q)pq^{-1}}$.  As $\pi$ is an
admissible action of $G$ on $B$, we have
$\pi_{qp^{-1}\rho_p(q^{-1})} = \pi_{p^{-1}}$. Now we easily obtain
that $((a \bowtie b)^\ast)^\ast = a \bowtie b$. \hfill
$\blacksquare$\\
\end{itemize}

Take the notations and assumptions as in Theorem 3.8.  Because
$D^\pi =  A^{cop} \bowtie \widetilde{B}$ is a twisted tensor
product multiplier Hopf algebra of $A^{cop}$ and $\widetilde{B}$,
it is quite obvious that integrals on $A$ and on $\widetilde{B}$
compose to an integral on $D^\pi = A^{cop} \bowtie \widetilde{B}$
in the following way.  Let $\varphi_A$ be an left integral on $A$
and let $\psi_B$ be a right integral on $B$. Consider
$\widetilde\psi$, defined on $B$ by
$\widetilde\psi(b)=\psi(\pi_{p^{-1}}(b))$ when $b\in B_p$. Then
$\varphi_A \otimes \widetilde{\psi}_B$ is a right integral on
$D^\pi = A^{cop} \bowtie \widetilde{B}$, see also Theorem 2.11(3).
\\ \ \\
We now consider the $^\ast$-situation. From [D, Remark 3.11] we
know that in this case, positive integrals on $A^{cop}$ and
$\widetilde{B}$ do not compose in a trivial way to a positive
integral on $D^\pi = A^{cop} \bowtie \widetilde{B}$.  In [De-VD,
Theorem 3.4], the problem for the usual Drinfel'd double $D =
A^{cop} \bowtie B$, which is associated to the multiplier Hopf
$^\ast$-algebra pairing $\langle A, B \rangle$, is treated as
follows. Let $\delta_A$ and $\delta_B$ denote the modular
multipliers in $M(A)$ and $M(B)$ respectively (see Theorem 1.5).
In [De-VD] there is given a meaning to the complex number $\langle
\delta_A, \delta_B\rangle^{1/2}$.  Furthermore, it is proven in
[De-VD, Theorem 3.4] that $\langle \delta_A, \delta_B\rangle^{1/2}
(\varphi_A \otimes \psi_B)$ is a positive right integral on $D =
A^{cop} \bowtie B$ whenever $\varphi_A$ is a positive left
integral on $A$ and $\psi_B$ is a positive right integral on $B$.
\\
In the following proposition we obtain this result for the
Drinfel'd double construction $D^\pi$.
\\ \ \\
{\bf 3.9 Proposition}\hs Let $A$ and $B$ be multiplier Hopf
$^\ast$-algebras as in Theorem 3.8(2).  Let $\varphi_A$ be a
positive left integral on $A$ and $\psi_B$ is a positive right
integral on $B$. Define as before $\widetilde\psi$ by
$\widetilde\psi(b)=\psi(\pi_{p^{-1}}(b))$ when $b\in B_p$. Then
$\langle \delta_A, \delta_B\rangle^{1/2} (\varphi_A \otimes
\widetilde{\psi}_B)$ is a positive right
integral on $D^\pi =  A^{cop} \bowtie \widetilde{B}$ .\\
\ \\
{\bf Proof.} \hs A straightforward calculation shows that
$$(\iota_A \otimes \widetilde{\psi}_B) (R(b \otimes a)) = \widetilde{\psi}_B(b)
(\delta_B^{-1} \btr a)$$ for all $a \in A$ and $b \in B$. Now, the
proof of [De-VD, Theorem 3.4] can be repeated, with the twist map
$R$  given by the formula $R(b \otimes a) = \sum (\pi_{p^{-1}}
(b_{(1)}) \btr a \btl S^{-1} (b_{(3)})) \otimes
b_{(2)}$ for all $ a \in A$ and $b \in B_p$.\insp $\blacksquare$\\
\ \\
{\it Two special cases}
\\ \ \\
The first case is the one with the trivial action. We get the
following (expected) result.
\\ \ \\
{\bf 3.10 Proposition} \hs Take the notations as in Theorem 3.8.
If the admissible action $\pi$ is the trivial action, then
$D^\pi$, constructed in Theorem 3.8, is nothing else but the usual
Drinfel'd double, $D = A^{cop} \bowtie B$, associated with the
pair $\langle A, B \rangle$ (as constructed and studied in
[Dr-VD] and [De-VD]).\\
\ \\
The other case is more interesting. Now, let $\pi$ an admissible
action such that for all $p, q \in G$ we have $\pi_p (B_q) =
B_{pqp^{-1}}$. This is the case with the adjoint action as in
Example 2.7. In the framework of Hopf group-coalgebras, these
actions are called crossings, see [T, Section 11]. We generalize
this definition here.\\
\ \\
{\bf  3.11 Definition}\hs Let $B$ be a regular $G$-cograded
multiplier Hopf algebra.  An admissible action $\pi$ of $G$ on $B$
is called a {\it crossing} if for all $p, q \in G$ we have
$\pi_p(B_q) =
B_{pqp^{-1}}$.\\
\ \\
In the following propositions, we describe the multiplier Hopf
algebras $\widetilde{B}$ and $D^\pi$ in more detail for the case
where the admissible action is a crossing.  We prove that the
considered multiplier Hopf algebras are again $G$-cograded.
Furthermore, we show that there is again a crossing of $G$ on
$\widetilde{B}$ and on $D^\pi$, defined in a natural way.
\\ \ \\
The following proposition generalizes the {\it mirror}
construction in the framework of crossed Hopf group-coalgebras,
see [T, Section 11].\\
\ \\
{\bf  3.12 Proposition}\hs Let $B$ be a regular $G$-cograded
multiplier Hopf algebra.  Let $\pi$ be a crossing of $G$ on $B$.
The deformed multiplier Hopf algebra $\widetilde{B} = (B,
\widetilde{\Delta})$ is again $G$-cograded.  Furthermore, $\pi$ is
also a crossing of $G$ on $\widetilde{B}$.  The deformation
$\widetilde{\widetilde{B}} = (B, \widetilde{\widetilde{\Delta}})$
equals the original $G$-cograded multiplier Hopf algebra $B$.\\
\ \\
{\bf Proof.} \hs From Theorem 2.11 we have that $\widetilde{B} =
(B, \widetilde{\Delta})$ is a regular multiplier Hopf algebra. Put
$\widetilde{B}_p= B_{p^{-1}}$. Then $\widetilde{B} =
\bigoplus\limits_{p\in G} \widetilde{B}_p$ and we have
\begin{eqnarray*}
\begin{array}{l}
\widetilde{\Delta} (\widetilde{B}_p) (1 \otimes \widetilde{B}_q) =
\widetilde{\Delta} (B_{p^{-1}}) (1 \otimes B_{q^{-1}}) =
(\pi_q\otimes \iota) (\Delta (B_{p^{-1}})(1 \otimes B_{q^{-1}}))\\
\hspace*{3cm}=(\pi_q\otimes \iota)(B_{p^{-1} q} \otimes
B_{q^{-1}}) = B_{qp^{-1}} \otimes B_{q^{-1}} =
\widetilde{B}_{pq^{-1}} \otimes \widetilde{B}_q.
\end{array}
\end{eqnarray*}
Similarly,
\begin{eqnarray*}
\begin{array}{l}
(\widetilde{B}_q\otimes 1) \widetilde{\Delta}(\widetilde{B}_p) =
(B_{q^{-1}} \otimes 1) \widetilde{\Delta} (B_{p^{-1}}) =
(B_{q^{-1}} \otimes 1) ((\pi_{q^{-1}p} \otimes \iota) \Delta
(B_{p^{-1}}))\\
\hspace*{2cm}= (\pi_{q^{-1}p} \otimes \iota) ((B_{p^{-1}q^{-1}p}
\otimes 1) \Delta (B_{p^{-1}})) = (\pi_{q^{-1}p} \otimes \iota)
(B_{p^{-1}q^{-1}p} \otimes B_{p^{-1}q})\\
\hspace*{2cm}= B_{q^{-1}} \otimes B_{p^{-1}q} =
\widetilde{B}_q\otimes \widetilde{B}_{q^{-1}p}.
\end{array}
\end{eqnarray*}
As the algebra structure of $\widetilde{B}$ is the same as the
algebra structure of $B$, $\pi$ is a crossing on $\widetilde{B}$
if for all $p \in G$ we have that $\pi_p$ respects the
comultiplication $\widetilde{\Delta}$ in the sense that
$\widetilde{\Delta}(\pi_p(b)) = (\pi_p \otimes \pi_p)
(\widetilde{\Delta}(b))$ for all $b \in B$.  To show this, take $b
\in B$ and $b' \in B_q$, then we have
\begin{eqnarray*}
\begin{array}{l}
\widetilde{\Delta}(\pi_p(b)) (1 \otimes b') = (\pi_{q^{-1}}
\otimes \iota) (\Delta(\pi_p(b)) (1 \otimes b')) = (\pi_{q^{-1}p}
\otimes \pi_p) (\Delta(b) (1 \otimes \pi_{p^{-1}} (b'))\\
\hspace*{1cm}= (\pi_p \otimes \pi_p) ((\pi_{p^{-1} q^{-1}p}
\otimes \iota) (\Delta(b) (1 \otimes \pi_{p^{-1}} (b')))) = (\pi_p
\otimes \pi_p)
(\widetilde{\Delta} (b) (1 \otimes \pi_{p^{-1}} (b')))\\
\hspace*{1cm}= ((\pi_p \otimes \pi_p) \widetilde{\Delta}(b)) (1
\otimes b').
\end{array}
\end{eqnarray*}
So $\widetilde{B} = (B, \widetilde{\Delta})$ is a $G$-cograded
multiplier Hopf algebra and $\pi$ is a crossing of $G$ on
$\widetilde{B}$.\\ Therefore we can consider the deformation
$\widetilde{\widetilde{B}} = (B, \widetilde{\widetilde{\Delta}})$.
We will now show that this is the original $G$-cograded multiplier
Hopf algebra $B$.  Take $b \in B$ and $b' \in \widetilde{B}_q =
B_{q^{-1}}$. Then we have
\begin{eqnarray*}
\widetilde{\widetilde{\Delta}} (b) (1 \otimes b') = (\pi_{q^{-1}}
\otimes \iota)(\widetilde{\Delta}(b) (1 \otimes b')) =
(\pi_{q^{-1}} \otimes \iota) ((\pi_q \otimes \iota) (\Delta(b) (1
\otimes b')) = \Delta(b) (1 \otimes b').
\end{eqnarray*}
Therefore we have $(B, \widetilde{\widetilde{\Delta}}) = B$.\insp
$\blacksquare$\\
\ \\
Next we consider the Drinfel'd double construction $D^\pi =
A^{cop} \bowtie \widetilde{B}$ when $\pi$ is a crossing.\\
\ \\
{\bf 3.13 Proposition}\hs Take the notations and assumptions as in
Theorem 3.8. Assume furthermore that the admissible action $\pi$
of $G$ on $B$ is a crossing.  Then $D^\pi =  A^{cop} \bowtie
\widetilde{B}$ is again a $G$-cograded multiplier Hopf algebra.
Furthermore, there is a crossing of $G$ on $D^\pi$, defined in a natural way.
\\
{\bf Proof.} \hs Let $D^\pi_p= A \bowtie B_{p^{-1}}$ for any $p\in
G$. Because $\pi$ is assumed to be a crossing, one can show, with
the techniques used before, that $R(B_p \otimes A)=A \otimes B_p$
for all $p\in G$. It follows easily that $D^\pi_p$ is a subalgebra
of $D^\pi$ and that $D^\pi_p D^\pi_q= 0$ if $p\neq q$. Also $D^\pi
= \bigoplus\limits_{p\in G} D^\pi_p$. This gives the desired
decomposition of the algebra $D^\pi$.
\\
Next we prove that $\overline{\Delta} (D^\pi_p) ((1 \bowtie 1)
\otimes D^\pi_q) = D^\pi_{pq^{-1}} \otimes D^\pi_q$.  Take $a, a'
\in A$ and $b \in B_{p^{-1}}$, $b' \in B_{q^{-1}}$.  Then, we have
\begin{eqnarray*}
\overline{\Delta} (a \bowtie b)((1 \bowtie 1) \otimes (a' \bowtie
b')) = \textstyle\sum (a_{(2)} \bowtie \pi_q(b_{(1)})) \otimes
((a_{(1)} \bowtie b_{(2)}) (a' \bowtie b')).
\end{eqnarray*}
As $\pi$ is a crossing of $G$ on $B$, we have that $b_{(2)} \in
B_{q^{-1}}$ and $\pi_q(b_{(1)}) \in B_{qp^{-1}}$. Similarly, to
prove that also $(D^\pi_p\otimes (1 \bowtie 1)) \overline{\Delta}
(D^\pi_q) = D^\pi_p\otimes D^\pi_{p^{-1}q}$, take $a$, $a'$, $b$
and $b'$ as above and write
\begin{eqnarray*}
((a \bowtie b) \otimes (1 \bowtie 1)) \overline{\Delta} (a'
\bowtie b') = \textstyle\sum ((a \bowtie b) (a'_{(2)} \bowtie
\pi_{p^{-1}q}(b'_{(1)}))) \otimes (a'_{(1)} \bowtie b'_{(2)})).
\end{eqnarray*}
Now we find that $b'_{(2)} \in B_{q^{-1}p}$ and $\pi_{p^{-1}q}
(b'_{(1)}) \in B_{p^{-1}}$. So, we have shown that $D^\pi$ is
$G$-cograded.
\\
Next, we will define a crossing of $G$ on $D^\pi$. First consider
the action $\pi'$ of $G$ on $A$ defined in the following way. Take
$p \in G$ and define the linear map $\pi'_p$ on $A$ by the formula
$\langle \pi'_p(a), b\rangle = \langle a, \pi_{p^{-1}} (b)\rangle$
for all $a \in A$ and $b \in B$. Clearly $\pi'_p$ is a linear
isomorphism such that $(\pi'_p)^{-1} = \pi'_{p^{-1}}$. From the
definition of $\pi'_p$, it easily follows that $\pi'_p$ is an
algebra isomorphism on $A$ and $\Delta(\pi'_p(a)) = (\pi'_p
\otimes \pi'_p) (\Delta(a))$. Furthermore we have $\pi'_{pq} =
\pi'_p \pi'_q$ for all $p, q \in G$.
\\
Next, for $p \in G$, consider the linear isomorphism $\pi'_p
\otimes \pi_p$ on $D^\pi = A^{cop} \bowtie \widetilde{B}$. It is
not difficult to see that $\pi'_p \otimes \pi_p$ is an algebra
isomorphism on $D^\pi$. Furthermore $\pi'_{pq} \otimes \pi_{pq} =
(\pi'_p \otimes \pi_p) (\pi'_q \otimes \pi_q)$ for all $p, q \in
G$.  Moreover $(\pi'_p \otimes \pi_p) (D^\pi_q) =
D^\pi_{pqp^{-1}}$ for all $p, q\in G$.\\
To complete the proof, we show that for all $p \in G$, the
isomorphism $\pi'_p \otimes \pi_p$ on $D^\pi$ respects the
comultiplication of $D^\pi$.  For any $a \in A$ and $b \in B$ we
have
\begin{eqnarray*}
\begin{array}{l}
\overline{\Delta} ((\pi'_p \otimes \pi_p) (a\bowtie b)) =
\Delta^{cop} (\pi'_p(a)) \widetilde{\Delta} (\pi_p(b)) \\
\hspace*{2cm} = (\pi'_p \otimes \pi'_p) (\Delta^{cop} (a)) (\pi_p
\otimes \pi_p)(\widetilde{\Delta} (b)) = ((\pi'_p \otimes \pi_p)
\otimes (\pi'_p \otimes \pi_p)) (\overline{\Delta} (a \bowtie b)).
\end{array}
\end{eqnarray*}
We conclude that the isomorphisms $\pi'_p \otimes \pi_p$, with $p
\in G$, define a crossing of $G$ on $D^\pi$. \hfill
$\blacksquare$\\
\ \\
{\bf 3.14 Example} \hs Let $G$ be any group.  Consider a Hopf
$G$-coalgebra as given in [T-Section 11].  In Proposition 2.2 we
saw that we can associate a regular $G$-cograded multiplier Hopf
algebra. We use the notations and assumptions of this proposition
and we put $B= \bigoplus\limits_{p\in G} B_p$. We now suppose that
each algebra $B_p$ is finite-dimensional. Next let $B^\ast =
\bigoplus\limits_{p\in G} (B_p)'$ where $(B_p)'$ the dual space of
$B_p$.  This is called the reduced dual space of $B$ and in
general, $B^\ast$ is smaller than the (full) dual of $B$. It can
be shown that $B^\ast$ is a Hopf algebra. The product is defined
dual to the coproduct as follows. Take $f \in (B_p)'$ and $g \in
(B_q)'$, then $fg$ in $(B_{pq})'$ is defined by the formula
$(fg)(x) = (f \otimes g) \Delta_{p,q}(x)$ for all $x \in B_{pq}$.
The unit in $B^\ast$ is $\varepsilon$ (on $B_1$). The coproduct is
defined dual to the product. For $f \in (B_p)'$ we have $\Delta
(f) \in (B_p)' \otimes (B_p)'$ when $\Delta$ is the dual to the
multiplication in $B_p$. For $f \in (B_p)'$ we have
$\varepsilon(f) = f(1_p)$ and $S(f) = f \circ S$.
\\
The evaluation map defines a pairing $\langle B^\ast, B\rangle$
between $B^\ast$ and $B$ of the type considerd in Proposition 3.1.
A crossing in the sense of [T - Section 11] gives a crossing of
$G$ on $B$ in the sense of Definition 3.11.  The Drinfel'd double
$D^\pi = (B^\ast)^{cop} \bowtie \widetilde{B}$ is a $G$-cograded
multiplier Hopf algebra, $D^\pi = \bigoplus\limits_{p\in G}
D^\pi_p$ with $D^\pi_p= B^\ast \bowtie B_{p^{-1}}$. For each $p
\in G$, we define the isomorphism $\pi'_p \in Aut(B^\ast)$ by the
formula $\langle \pi'_p(f), b \rangle = \langle f, \pi_{p^{-1}}
(b)\rangle$.  From Proposition 3.12 we have that for all $p \in
G$, the isomorphisms $\pi'_p \otimes \pi_p$ on $D^\pi =
(B^\ast)^{cop} \bowtie \widetilde{B}$ provide a crossing of the
group $G$ on $D^\pi$.\\ This Drinfel'd double $D^\pi$ is the same
as the one constructed in [Z, Section 5] in the framework of Hopf
group-coalgebras.


\begin{thebibliography}{AAAAAA}
\bibitem[A-D-VD]] A.T.\ Abd El-Hafez, L.\ Delvaux and A.\ Van
Daele, {\it Group-cograded multiplier Hopf ($^*$-)algebras},
Preprint (2004). Sever version math.QA/04????.
\bibitem[D]] L.\ Delvaux, {\it Twisted tensor product of multiplier
Hopf ($^\ast$-)algebras}, J.\ Algebra {\bf 269}(1) (2003),
285--316.
\bibitem[De-VD]] L.\ Delvaux and A. Van Daele, {\it The Drinfel'd
double of multiplier Hopf algebras}, J.\ Algebra {\bf 272} (1)
(2004), 273--291.
\bibitem[Dr-VD]] B.\ Drabant, A.\ Van Daele, {\it Pairing and
quantum double of multiplier Hopf algebras}, Algebras and
Representative Theory {\bf 4}(2) (2001), 109--132.
\bibitem[D-VD-Z]] B.\ Drabant, A.\ Van Daele and Y.\ Zhang, {\it
Actions of multiplier Hopf algebras}, Commun. Algebra {\bf 27}(9)
(1999), 4117--4127.
\bibitem[H-A-M]] A.S.\ Hegazi, A.T.\ Abt El-Hafez and M.\ Mansour,
{\it Multiplier Hopf group-coalgebras}, Preprint Mansoura
University (2002).
\bibitem[K-V]] J.\ Kustermans and S.\ Vaes, {\it Locally compact
quantum groups}, Ann.\ Sci.\ Ec.\ Norm.\ Sup.\ {\bf 33} (2000),
837--934.
\bibitem[K-VD]] J.\ Kustermans and A.\ Van Daele, {\it
C$^*$-algebraic quantum groups arising from algebraic quantum
groups}, Int.\ J.\ Math.\ {\bf 8} (1997), 1067--1139.
\bibitem[T]] V.G.\ Turaev, {\it Homotopy field theory in
dimension 3 and crossed group-categories}, Preprint GT/0005291.
\bibitem[VD1]] A.\ Van Daele, {\it Multiplier Hopf algebras},
Trans.\ Am.\ Math.\ Soc.\ {\bf 342}(2) (1994), 917--932.
\bibitem[VD2]] A.\ Van Daele, {\it An algebraic framework for
group duality}, Advances in  Ma\-thematics {\bf 140} (1998),
323--366.
\bibitem[VD3]] A. Van Daele, {\it Multiplier Hopf $^\ast$-algebras
with positive integrals: A laboratory for locally compact quantum
group}.  In ``Locally Compact Quantum Groups and Groupoids'',
Proceedings of the Meeting of Theoretical Physicists and
Mathematicians, Strasbourg, February 21-23, 2003 (editor L.
Vainerman).
\bibitem[VD-Z]] A.\ Van Daele and Y.\ Zhang, {\it A survey on
multiplier Hopf algebras}, in: S.\ Caenepeel, F.\ Van Oystayen
(Eds.), Hopf Algebras and Quantum Groups, Dekker, New York, 1998
pp.259-309.
\bibitem[Z]] M.\ Zunino, {\it Double construction for
crossed Hopf coalgebras}, Preprint QA/0212192.
\end{thebibliography}
\end{document}